\documentclass{article}

\usepackage{arxiv}
\usepackage{booktabs}
\usepackage{siunitx}
\usepackage[utf8]{inputenc} 
\usepackage[T1]{fontenc}    
\usepackage{hyperref}       
\usepackage{url} 
\usepackage{amsmath}
\usepackage{cite}
\usepackage{booktabs}       
\usepackage{amsfonts}       
\usepackage{nicefrac}       
\usepackage{microtype}      
\usepackage{lipsum}
\usepackage{graphicx}
\usepackage{overpic}

\title{Safe cross-entropy-based importance sampling for rare event simulations}

\date{}

\author{
Zhiwei Gao\thanks{First author. Email: \texttt{zhiwei\_gao@brown.edu}} \\
Applied Mathematics \\
Brown University \\
\texttt{zhiwei\_gao@brown.edu}
\And
George Karniadakis\thanks{Corresponding author. Email: \texttt{george\_karniadakis@brown.edu}} \\
Applied Mathematics \\
Brown University \\
\texttt{george\_karniadakis@brown.edu}
}

\begin{document}
\maketitle
\begin{abstract}
The Improved Cross‐Entropy (ICE) method is a powerful tool for estimating failure probabilities in reliability analysis. Its core idea is to approximate the optimal importance‐sampling density by minimizing the forward Kullback–Leibler divergence within a chosen parametric family—typically a mixture model. However, conventional mixtures are often light‐tailed, which leads to slow convergence and instability when targeting very small failure probabilities. Moreover, selecting the number of mixture components in advance can be difficult and may undermine stability. To overcome these challenges, we adopt a weighted cross‐entropy–penalized expectation–maximization (EM) algorithm that automatically prunes redundant components during the iterative process, making the approach more stable. Furthermore, we introduce a novel two‐component mixture that pairs a light‐tailed distribution with a heavy‐tailed one, enabling more effective exploration of the tail region and thus accelerating convergence for extremely small failure probabilities.  We call the resulting method Safe-ICE and assess it on a variety of test problems. Numerical results show that Safe-ICE not only converges more rapidly and yields more accurate failure‐probability estimates than standard ICE, but also identifies the appropriate number of mixture components without manual tuning.   
\end{abstract}


\section{Introduction}
Reliability analysis \cite{gardoni2017risk,zacks2012introduction} aims to quantify the probability that an engineering system or component enters a predefined failure state. Mathematically, this is often formulated in terms of a random input vector $\mathbf{u} \in \Omega \subseteq \mathbb{R}^d$, representing uncertain system parameters with a prior distribution $p(\mathbf{u})$. A failure event is typically characterized by a \textit{limit-state function} (LSF) $g(\mathbf{u})$, with failure occurring when $g(\mathbf{u}) \leq 0$. That is, the failure domain is defined as
$
\Omega_{\mathcal{F}} = \{\mathbf{u} \in \Omega : g(\mathbf{u}) \leq 0\}.
$
The function $g$ may represent one or more distinct failure modes, corresponding to either component-level or system-level reliability considerations. The probability of failure is then expressed as a $d$-dimensional integral over the failure region:
\begin{equation}
\label{failure_probability}
P_\mathcal{F} = \mathbb{P}(\Omega_{\mathcal{F}}) = \int_{\Omega} \mathbb{I}_{\Omega_{\mathcal{F}}}(\mathbf{u})\, p(\mathbf{u}) \, d\mathbf{u},
\end{equation}
where $\mathbb{I}_{\Omega_\mathcal{F}}:\Omega\to\{0,1\}$ represents the indicator function, which takes value $\mathbb{I}_{\Omega_\mathcal{F}}(\mathbf{u})=1$ when $\mathbf{u}\in\Omega_\mathcal{F}$ and $\mathbb{I}_{\Omega_\mathcal{F}}(\mathbf{u})=0$ otherwise. The prior is typically Gaussian; otherwise, a Nataf \cite{lebrun2009innovating} or Rosenblatt \cite{lebrun2009rosenblatt} transformation can be applied to map the original distributions to Gaussian ones.

Analytically computing the failure probability is generally infeasible because it involves a high-dimensional integral, and the limit-state function 
$g$ is computationally intensive when it is governed by partial differential equations (PDEs). To alleviate these issues, different methods have been developed to reduce the number of LSF evaluations and increase the accuracy as well as stability. These methods generally can be divided into three categories: i) optimization approaches, such as First/Second order reliability analysis \cite{maier2001first,kiureghian1991efficient}, ii) surrogate models \cite{bichon2011efficient,li2010evaluation,wei2019structural,ehre2022sequential,dubourg2011adaptive}, which replace $g$ with a less expensive model, and iii) Monte Carlo-based (MC) sampling methods \cite{cerou2012sequential,au2001estimation,homem2002estimation,denny2001introduction,kim2000nonparametric,papaioannou2015mcmc}. While optimization methods can be effective when $g$ is sufficiently smooth, sampling methods are more adaptive with complex geometries and exhibit better asymptotic convergence properties \cite{bucklew2004introduction}. 

Utilizing vanilla Monte Carlo to estimate the failure probability requires generating a set of samples from the prior $p(x)$ and then calculating the sample mean of the indicator function. However, MC simulation requires $\mathcal{O}(\frac{1}{\mathcal{P}_{\mathcal{F}}})$ to achieve a predefined coefficient of variation (CV), which can still be quite computationally demanding for very small failure probabilities. To alleviate this issue, importance sampling \cite{tokdar2010importance} (IS) provides an alternate way to reduce the variance by employing another biased distribution $q(\mathbf{u})$. To this end, Eq. \eqref{failure_probability} can be rewritten as 
\begin{equation}
\label{importance_failure_probability}
    P_\mathcal{F} = \int_{\Omega} \mathbb{I}_{\Omega_{\mathcal{F}}}(\mathbf{u})\, p(\mathbf{u}) \, d\mathbf{u} = \mathbb{E}_{q}[\mathbb{I}_{\Omega_{\mathcal{F}}}(\mathbf{u})R(\mathbf{u})],
\end{equation}
where $R(\mathbf{u}) = \frac{p(\mathbf{u})}{q(\mathbf{u})}$ is the weight function.  Subsequently, by generating a set of samples $\{\mathbf{u}_{i}\}_{i=1}^{N}$ from $q(\mathbf{u})$, one can obtain an unbiased estimate $\widehat{P}_{\mathcal{F}}$ as 
\begin{equation*}
    \widehat{P}_{\mathcal{F}} = \frac{1}{N}\sum_{i=1}^{N}\mathbb{I}_{\Omega_{\mathcal{F}}}(\mathbf{u}_{i})R(\mathbf{u}_{i}),
\end{equation*}
provided that $\text{supp}(p) \subseteq \text{supp}(q)$. Using the variational method, the theoretically optimal $p^{*}(\mathbf{u})$ with zero variance can be derived as 
\begin{equation}
    \label{optimal_bias_density}
    p^{*}(\mathbf{u}) = \frac{\mathbb{I}_{\Omega_{\mathcal{F}}}(\mathbf{u})p(\mathbf{u})}{P_{\mathcal{F}}}.
\end{equation}

However, the optimal IS density $p^{*}(\mathbf{u})$ is not available in practice due to the unknown $P_{\mathcal{F}}$. Therefore, the common approach is to find an approximately optimal $q^{*}(\mathbf{u})$ in a sequential manner \cite{papaioannou2016sequential,beaurepaire2013reliability},  as in the adaptive importance sampling method. In such cases, the biased density $q(\mathbf{u})$ is frequently chosen as a non-parametric \cite{zhang1996nonparametric} or mixture parametric distribution \cite{cappe2008adaptive} to adapt to the unimodal or multimodal failure region.

The cross-entropy (CE) \cite{de2005tutorial,kurtz2013cross,papaioannou2019improved} method is a principled, adaptive importance-sampling framework that iteratively refines a mixture parametric proposal by minimizing the forward KL divergence over a sequence of increasingly rare intermediate failure regions. When the proposal belongs to an exponential family, these updates typically admit closed-form Expectation-Maximization (EM) steps \cite{teng2010algorithm}. However, CE discards most samples by relying solely on an elite subset—diminishing statistical efficiency—and, in high dimensions, its Gaussian proposals collapse onto a thin “shell” \cite{ledoux2001concentration},  causing numerical instability. The improved CE (ICE) method remedies this by replacing the hard indicator with a smoothed surrogate to leverage all samples and by employing a von Mises–Fisher–Nakagami (vMFNM) mixture \cite{papaioannou2019improved} for greater tail flexibility. However, ICE still requires a user-specified number of mixture components and lacks an intrinsic mechanism for adapting this choice. Furthermore, while ICE converges quickly for moderately small failure probabilities, its performance degrades when probabilities become vanishingly small due to the light‐tailed Nakagami kernels. This results in a substantial increase in the number of iterations, thereby raising the overall computational burden. Moreover, the estimation becomes more unstable, particularly in cases where the failure region exhibits multi‐modal features.

To address the aforementioned challenges, we adopt a weighted cross-entropy–penalized EM algorithm, where the additional regularization in the M-step automatically prunes redundant mixture components and stabilizes the update process~\cite{yang2012robust}. Building on the idea of safe adaptive importance sampling~\cite{delyon2021safe}, we further design a mixture proposal that combines the vMFNM density~\cite{papaioannou2019improved} with a dedicated heavy-tailed distribution. This construction accelerates convergence while maintaining robustness for extremely small failure probabilities. At each iteration, samples are drawn from the mixture, and the vMFNM parameters are updated with the penalized EM algorithm, progressively reducing components and driving the proposal more efficiently and stably toward the optimal biasing distribution. To balance exploration and exploitation, the weight of the heavy-tailed distribution is gradually decreased following a cosine-annealing schedule. Extensive experiments demonstrate that the proposed method significantly improves accuracy and stability while reducing computational cost compared to the baseline ICE algorithm. Our contributions can be summarized as follows:
\begin{itemize}
     \item We adopt a weighted cross‐entropy–penalized EM algorithm to remove redundant components from the mixture model, and derive a new weight‐update rule that promotes sparsity during the M-step.
    \item We propose a novel biased distribution that integrates both light‐tailed and heavy‐tailed components to accelerate convergence. Furthermore, the proportion of the heavy‐tailed component is adaptively adjusted via a cosine annealing strategy, leading to consistently strong performance.
    \item To enhance the quality of generated samples, we adopt the inverse Nakagami distribution in combination with the von Mises-Fisher distribution as the heavy‐tailed component.
    \item We validate the effectiveness of the proposed techniques through a series of numerical experiments. The results show that our method achieves higher accuracy and greater stability than the original ICE approach, while also reducing computational cost.
    
\end{itemize}

The remainder of this paper is organized as follows. Section 2 introduces the improved cross-entropy framework for rare event simulation. Section 3 introduces our novel mixtures as well as the cross-entropy penalized EM. To verify the effectiveness of our proposed framework, we implement different experiments in Section 4.

\section{Improved Cross-Entropy Method}
The CE method \cite{kurtz2013cross} and its improved version (ICE) \cite{papaioannou2019improved} are adaptive importance sampling schemes that aim to approximate the optimal IS density in \eqref{optimal_bias_density} within a specified parametric family. Let $q(\mathbf{u};\mathbf{v})$ be the predefined parametric family with parameters $\mathbf{v}$. CE (ICE) aims to find the optimal parameters $\mathbf{v}^{*}$ by minimizing the following KL divergence,
\begin{equation}
\label{KL_divergence}
    \mathbb{D}[p^{*}\|q] = \mathbb{E}_{p^{*}}[\ln p^{*}(\mathbf{u})] - \mathbb{E}_{p^{*}}[\ln q(\mathbf{u};\mathbf{v})].
\end{equation}
By discarding the first term and substituting Eq. \eqref{optimal_bias_density} into Eq. \eqref{KL_divergence}, the optimization is equivalent to finding the optimal parameters via the following maximization problem:
\begin{equation}
\label{maximation_problem}
\begin{split}
    \mathbf{v}^{*} &= \arg\max_{\mathbf{v}\in V} \mathbb{E}_{p^{*}}[\ln q(\mathbf{u};\mathbf{v})]\\ 
    & = \arg\max_{\mathbf{v}\in V}\int_{\Omega}\mathbb{I}_{\Omega_{\mathcal{F}}}(\mathbf{u})p(\mathbf{u})\ln q(\mathbf{u};\mathbf{v})d\mathbf{u}, 
\end{split}
\end{equation}
where $V$ is the parameter space. Solving problem \eqref{maximation_problem} requires approximating the high-dimensional integral using the MC simulation. In detail, suppose $\{\mathbf{u}_{i}\}_{i=1}^{N}\sim p(\mathbf{u})$, Eq. \eqref{maximation_problem} can be approximated as 
\begin{equation}
\label{discrete_optimization}
    \hat{\mathbf{v}}^{*} = \arg\max_{\mathbf{v}\in V}\sum_{i=1}^{N}\mathbb{I}_{\Omega_{\mathcal{F}}}(\mathbf{u}_{i})\ln q(\mathbf{u}_{i};\mathbf{v}).
\end{equation}
However, the MC estimate could lead to a large variance. Instead, importance sampling chooses another biased density $q(\mathbf{u};\mathbf{w})$ from the same parametric family and generate $\{\mathbf{u}_{i}\}_{i=1}^{N}\sim q(\mathbf{u};\mathbf{w})$, leading to the following optimization problem,
\begin{equation}
\label{importance_discrete_optimization}
    \hat{\mathbf{v}}^{*} = \arg\max_{\mathbf{v}\in V}\sum_{i=1}^{N}\mathbb{I}_{\Omega_{\mathcal{F}}}(\mathbf{u}_{i})\ln q(\mathbf{u}_{i};\mathbf{v})\frac{p(\mathbf{u}_{i})}{q(\mathbf{u}_{i};\mathbf{w})}.
\end{equation}
Nonetheless, accurately evaluating the integral in \eqref{maximation_problem} requires a significant number of samples that fall in $\Omega_{\mathcal{F}}$, which is infeasible when $P_{\mathcal{F}}$ is critically small. One possible solution is to introduce a sequence of intermediate distributions $\{p_{t}\}_{t=1}^{T}$ to gradually transit the prior $p(\mathbf{u})$ to the IS optimal density $p^{*}$. Hence, the optimization problem \eqref{importance_discrete_optimization} can be decomposed as a sequence of intermediate optimization problems, i.e., 
\begin{equation}
\label{intermediate_optimization_problems}
\hat{\mathbf{v}}_{t} = \arg\max_{\mathbf{v}\in V}\sum_{i=1}^{N}\ln q(\mathbf{u}_{i};\mathbf{v})W_{t}(\mathbf{u}_{i}, \hat{\mathbf{v}}_{t-1}),
\end{equation}
where $\{\mathbf{u}\}_{i=1}^{N} \sim q(\mathbf{u};\widehat{\mathbf{v}}_{t-1})$ and $W_{t}(\mathbf{u}_{i}, \widehat{\mathbf{v}}_{t-1}) = \frac{p_{t}(\mathbf{u}_{i})}{q(\mathbf{u}_{i};\mathbf{v}_{t-1})}$ is the weight function between the intermediate distribution and the last bias distribution.
Under the setting of ICE, the intermediate distribution $p_{t}(\mathbf{u})$ is chosen as 
\begin{equation}
\label{intermediate_diatribution}
p_{t}(\mathbf{u}) = \frac{1}{P_{t}}\Phi(-\frac{g(\mathbf{u})}{\sigma_{t}})p(\mathbf{u}) = \frac{1}{P_{t}}h_{t}(\mathbf{u})p(\mathbf{u}),
\end{equation}
where $\Phi(\cdot)$ is the cumulative density function of the normal distribution used to approximate the indicator function, $\sigma_{t}$ is a sequence of decreasing smoothing parameters, and $P_{t}$ is the intermediate failure probability. In detail, we have 
\begin{equation*}
    \lim_{\sigma\to 0}\Phi(-\frac{g(\mathbf{u})}{\sigma}) = \mathbb{I}_{\Omega_{\mathcal{F}}}(\mathbf{u}). 
\end{equation*}
Therefore, as $\sigma_{t}\to 0$, $p_{t}$ converges to $p^{*}$. To this end, the ICE method can utilize all of the samples in each iteration to estimate the parameters of the next bias density, whereas the CE method only uses the elite samples, leading to wasted information. 

To ensure the smooth transition between $p_{t-1}(\mathbf{u})$ and $p_{t}(\mathbf{u})$, the smoothing parameter $\sigma_{t}$ at the next iteration can be adaptively chosen by the following minimization problem:
\begin{equation}
\label{sigma_optimization}
    \sigma_{t} = \arg\min_{\sigma \in (0, \sigma_{t-1})}(\delta_{W_{t}}(\sigma) - \delta_{\text{target}})^{2},
\end{equation}
where $\delta_{W_{t}}$ is the coefficient of variation (CV) of the weights $\{W_{t}(\mathbf{u}_{i})\}_{i=1}^{N}$ and $\delta_{\text{target}}$ is a predefined target value. 

In terms of the stopping criteria, ICE will terminate when the intermediate distribution $p_{t}$ is sufficiently close to $p^{*}$. However, since the KL divergence cannot be calculated analytically, an alternative choice using the importance weights 
\begin{equation}
    \label{stop_criteria}
    W^{*}_{t}(\mathbf{u}) = \frac{p^{*}(\mathbf{u})}{p_{t}(\mathbf{u})} \propto \frac{\mathbb{I}_{\Omega_{\mathcal{F}}}(\mathbf{u})}{h_{t}(\mathbf{u})}
\end{equation}
is applied. In detail, after using the samples from $q_{t}$ to calculate $\{W_{t}^{*}(\mathbf{u}_{i})\}_{i=1}^{N}$, the whole procedure of ICE can be stopped when the CV of $W^{*}_{t}$ is below some specified tolerance $\delta^{*}$. Generally, the whole procedure can be summarized as follows:
\begin{itemize}
    \item[1.] Initialize $\hat{\mathbf{v}}_{0}$ based on the specified parametric family. Initialize $\sigma_{0}$ and set $t = 0$. 
    \item[2.] Generate samples $\{\mathbf{u}_{i}\}_{i=1}^{N}\sim q(\mathbf{u};\hat{\mathbf{v}}_{t})$ and calculate the LSF values $g(\mathbf{u}_{i}), i = 1, \ldots, N$. 
    \item[3.] Calculate the weights $\{W_{t}^{*}(\mathbf{u}_{i})\}_{i=1}^{N}$ and the corresponding CV $\delta_{W_{t}^{*}}$. If it is larger than $\delta^{*}$, go to step 4. Otherwise, go to step 6.
    \item[4.] Solve optimization problem \eqref{sigma_optimization} to determine $\sigma_{t}$. 
    \item [5.] Solve the stochastic problem \eqref{importance_discrete_optimization} to determine the paramters $\hat{\mathbf{v}}_{t+1}$. Set $t = t+1$ and go to step 2.
    \item[6.] Approximate the failure probability $\widehat{P}_{\mathcal{F}}$ using the samples from the final stage. Set $T = t$, then 
    \begin{equation}
    \label{failure_estimation}
        \widehat{P}_{\mathcal{F}} = \frac{1}{N}\sum_{i=1}^{N}\mathbb{I}_{\Omega_{\mathcal{F}}}(\mathbf{u}_{i})\frac{p(\mathbf{u}_{i})}{q(\mathbf{u}_{i};\hat{\mathbf{v}}_{T})}.
    \end{equation}
\end{itemize}

For the choice of the parametric family, \cite{papaioannou2019improved} proposed a novel and flexible mixture model, the vMFNM, to approximate the optimal IS density $p^{*}$ instead of using simple Gaussian mixtures. This model generally performs well, even in high-dimensional settings. However, when dealing with small failure probabilities, the convergence rate of such light-tailed distributions can be slow, resulting in high computational costs for expensive LSF evaluations. Furthermore, if the failure region exhibits multi-modal features, the approximation may become highly unstable, possibly due to the weight degeneracy, especially when there is no prior knowledge of the number of modes. Therefore, inspired by sparse EM techniques and safe adaptive importance sampling, we propose a new mixture model to accelerate convergence. In particular, the EM framework is employed with a sparsity-inducing mechanism to automatically remove redundant modes during the updating process, thereby enhancing stability and improving the quality of the approximation.

\section{Safe Cross-Entropy Based Importance Sampling}
In this section, we will introduce the cross-entropy penalized EM to automatically erase the additional components, and we will present our new mixture model based on a heavy-tailed distribution. 
\subsection{Cross-entropy penalized EM}

To solve the stochastic optimization problem \eqref{intermediate_optimization_problems}, one needs to specify the parametric family $q$ first. Typically, $q$ is chosen to be a mixture model to match the multi-modal features of the failure region. In detail, suppose $q$ has the following form:
\begin{equation}
\label{mixture_family}
    q(\mathbf{u};\mathbf{v}) = \sum_{k=1}^{K}\pi_{k}q(\mathbf{u};\mathbf{v}_{k}), 
\end{equation}
where $\pi_{k}\geq 0$ are weights such that $\sum_{k=1}^{K}\pi_{k} = 1$, and $\mathbf{v} = \{\pi_{k}, \mathbf{v}_{k}\}_{k=1}^{K}$ is the complete set of parameters needed to be estimated. Suppose there is a set of samples $\mathcal{U} = \{\mathbf{u}_{i}\}_{i=1}^{N}\sim q(\mathbf{u};\widehat{\mathbf{v}}_{t-1})$, the weighted log likelihood in Eq. \eqref{intermediate_optimization_problems} can be written as 
\begin{equation}
    \label{incomplete_likelihood}
    l(\mathbf{v};\mathcal{U}) = \sum_{i=1}^{N}\log q(\mathbf{u}_{i};\mathbf{v})W_{t}(\mathbf{u}_{i};\widehat{\mathbf{v}}_{t-1}) = \sum_{i=1}^{N}W_{t}(\mathbf{u}_{i};\widehat{\mathbf{v}}_{t-1})\log \sum_{k=1}^{K} \pi_{k}q(\mathbf{u}_{i};\mathbf{v}_{k}).
\end{equation}
Generally, the maximum likelihood (ML) estimate by this incomplete likelihood cannot be obtained \cite{yang2012robust}. Instead, the EM algorithm employs an iterative procedure to update the parameters. In detail, the missing part is a set of labels $\mathcal{Z} = \{\mathbf{z}^{(1)}, \ldots, \mathbf{z}^{(N)}\}$, where each label is a binary vector $\mathbf{z}^{(i)} = [z_{1}^{(i)}, \ldots, z_{K}^{(i)}]$ with $z_{k}^{(i)} = 1$ and $z_{p}^{(i)} = 0$ for $p\neq k$, indicating the $\mathbf{u}_{i}$ is produced from the $k$th component. Then, the complete weighted log likelihood can be written as 
\begin{equation}
    \label{complete_likelihood}
    l(\mathbf{v};\mathcal{U}, \mathcal{Z}) = \sum_{i=1}^{N}W_{t}(\mathbf{u}_{i};\widehat{\mathbf{v}}_{t-1})\sum_{k=1}^{K}z_{k}^{(i)}\log\pi_{k}q(\mathbf{u}_{i};\mathbf{v}_{k}). 
\end{equation}
Therefore, EM can produce a sequence of estimates $\widehat{\mathbf{v}}(j)$ by alternatively applying two steps (until convergence is achieved): 
\begin{itemize}
    \item[1.] \textbf{E-step}: given the current estimate $\widehat{\mathbf{v}}(j)$, we have 
    \begin{equation}
    \begin{split}
        \label{expectation}
        \gamma_{k}^{(i)} &= \mathbb{E}[z_{k}^{(i)}|\mathcal{X}, \widehat{\mathbf{v}}(j)] = P(z_{k}^{(i)}=1|\mathbf{u}_{i}, \widehat{\mathbf{v}}(j))\\ 
        & = \frac{\widehat{\pi}_{k}(j)q(\mathbf{u}_{i};\widehat{\mathbf{v}}_{k}(j))}{\sum_{k=1}^{K}\widehat{\pi}_{k}(j)q(\mathbf{u}_{i};\widehat{\mathbf{v}}_{k}(j))}.
    \end{split}
    \end{equation}
    \item[2.] \textbf{M-step}: Update the parameters according to 
    \begin{equation}
        \label{m-step}
        \widehat{\mathbf{v}}(j+1) = \arg\max_{v\in V}Q(\mathbf{v}, \widehat{\mathbf{v}}(j)),
    \end{equation}
    where the so-called $Q$ function is defined as 
    \begin{equation}
    \begin{split}
        \label{Q-function}
        Q(\mathbf{v}, \widehat{\mathbf{v}}(j)) &= \mathbb{E}[l(\mathbf{v};\mathcal{X},\mathcal{Z})|\mathcal{X},\widehat{\mathbf{v}}(j)]\\ 
        & = \sum_{i=1}^{N}W_{t}(\mathbf{u}_{i};\widehat{\mathbf{v}}_{t-1})\sum_{k=1}^{K}\gamma_{k}^{(i)}\log\pi_{k}q(\mathbf{u}_{i};\mathbf{v}_{k}).
    \end{split}
    \end{equation}
\end{itemize}
By maximizing the $Q$ function, we can get the update formula of parameters $\pi_{k}, \mathbf{v}_{k}$. If $q$ is chosen to be the vMFNM distribution in \cite{papaioannou2019improved}, the M-step admits closed-form updates. In particular, the update for the weights can be written as 
\begin{equation}
\label{EM_step}
    \widehat{\mathbf{\pi}}_{k}^{\mathrm{EM}}(j+1) = \frac{\sum_{i=1}^{N}\gamma_{k}^{(i)}W_{t}(\mathbf{u}_{i};\widehat{\mathbf{v}}_{t-1})}{\sum_{i=1}^{N}\sum_{s=1}^{K}\gamma_{s}^{(i)}W_{t}(\mathbf{u}_{i};\widehat{\mathbf{v}}_{t-1})}.
\end{equation}
The whole procedure will stop until $|l(\widehat{\mathbf{v}}(j+1)) - l(\widehat{\mathbf{v}}(j))| \leq \epsilon \cdot l(\widehat{\mathbf{v}}(j))$, where $\epsilon$ is a predefined tolerance. After converging, we can set $\widehat{\mathbf{v}}_{t} = \widehat{\mathbf{v}}(j+1)$ to be the parameters of the next bias density. Note that the common EM algorithm cannot erase redundant mixture components, which may lead to numerical instability. To achieve this goal, inspired by \cite{yang2012robust}, we add a weighted cross-entropy term to the $Q$ function in the M-step to penalize the weights. Similarly, the weighted $Q$ function can be rewritten as 
\begin{equation}
    \label{novel_Qfunction}
    Q(\mathbf{v}, \widehat{\mathbf{v}}(j)) = \sum_{i=1}^{N}W_{t}(\mathbf{u}_{i};\widehat{\mathbf{v}}_{t-1})\sum_{k=1}^{K_{j}}\gamma_{k}^{(i)}\log\pi_{k}q(\mathbf{u}_{i};\mathbf{v}_{k}) + \beta\sum_{i=1}^{N}W_{t}(\mathbf{u}_{i};\widehat{\mathbf{v}}_{t-1})\sum_{k=1}^{K_{j}}\pi_{k}\ln \pi_{k}, \beta > 0.
\end{equation}
The update of $\pi_{k}$ can be derived (see appendix) by maximizing $Q$ under the constraint $\sum_{k=1}^{K_{j}}\pi_{k} = 1$, with the form of 
\begin{equation}
    \label{update_pi}
    \widehat{\pi}_{k}(j+1) = \widehat{\pi}_{k}^{\mathrm{EM}}(j+1) + \beta \left(\frac{\sum_{i=1}^{N}W_{t}(\mathbf{u}_{i};\widehat{\mathbf{v}}_{t-1})}{\sum_{i=1}^{N}\sum_{s=1}^{K_{j}}\gamma_{s}^{(i)}W_{t}(\mathbf{u}_{i};\widehat{\mathbf{v}}_{t-1})}\right)\widehat{\mathbf{\pi}}_{k}(j)\left(\ln \widehat{\mathbf{\pi}}_{k}(j) - \sum_{s=1}^{K_{j}}\widehat{\mathbf{\pi}}_{s}(j)\ln \widehat{\mathbf{\pi}}_{s}(j)\right).
\end{equation}
Using the update rule \eqref{update_pi}, some components will gradually decrease. Therefore, by discarding the components where $\widehat{\pi}_{k}(j+1) \leq 0$ according to \cite{yang2012robust}, we can update the number of components with 
\begin{equation}
\label{update_k}
    K_{j+1} = K_{j} - |\widehat{\pi}_{k}(j+1)|\widehat{\pi}_{k}(j+1) \leq 0, k = 1,\ldots, K_{j}|.
\end{equation}
Then, we can normalize $\widehat{\pi}_{k}(j+1)$ and $\gamma_{k}^{(i)}$ to continue updating the rest parameters using the M-step \eqref{m-step}. For the choice of $\beta$, to maintain a positive number of components, it is recommended to use the following formula based on \cite{yang2012robust}:
\begin{equation}
    \label{beta_update}
    \beta_{j+1} = \min\Biggl\{
\; \frac{1}{K_{j}}\sum_{k=1}^{K_{j}}
\exp\bigl(-\eta N\lvert \widehat{\pi}_{k}(j+1)
- \widehat{\pi}_{k}(j)\rvert\bigr)
\;,\;
\frac{\left(1 - \widehat{\pi}_{(1)}^{\mathrm{EM}}\right)}
{\left(-\,\widehat{\pi}_{(1)}(j)\,E\right)}
\Biggr\},
\end{equation}
where $\eta = \min\!\Bigl\{1,\;0.5^{\lfloor\tfrac{d}{2}-1\rfloor}\Bigr\}$, and 
\begin{equation} 
\widehat{\pi}_{(1)}^{\mathrm{EM}} = \max_{1\leq k\leq K_{j} }\widehat{\pi}_{k}^{\mathrm{EM}}(j+1), \,\widehat{\pi}_{(1)}(j) = \max_{1\leq k\leq K_{j}} \widehat{\pi}_{k}(j),\, E = \sum_{k=1}^{K_{j}}\widehat{\pi}_{k}(j)\ln \widehat{\pi}_{k}(j). 
\end{equation}
Note that larger values of $\beta$ impose a stronger penalty, which in turn drives the removal of more mixture components. In our ICE algorithm, we therefore reset $\beta = 1$ at the start of each outer loop and then update it according to Eq.~\eqref{beta_update}. For the selection of the initial $K$, we are supposed to choose a comparatively large $K$ such that an adequate estimate can still be obtained even with many ICE iterations. In practice, we found $K=20$ is enough to obtain good performance.

The choice of the parametric family can be diversified, but mainly focuses on the exponential family, including the mixture Gaussian distribution and the vMFNM distribution used in ICE. However, as these models are typically light-tailed, their convergence can be slow for extremely small failure probabilities, often requiring many ICE iterations. Therefore, we will introduce a novel mixture in the following section to accelerate the convergence. 

\subsection{Adaptive mixture with a heavy-tailed distribution}
\begin{figure}[htbp]
    \centering
\includegraphics[width=0.5\linewidth]{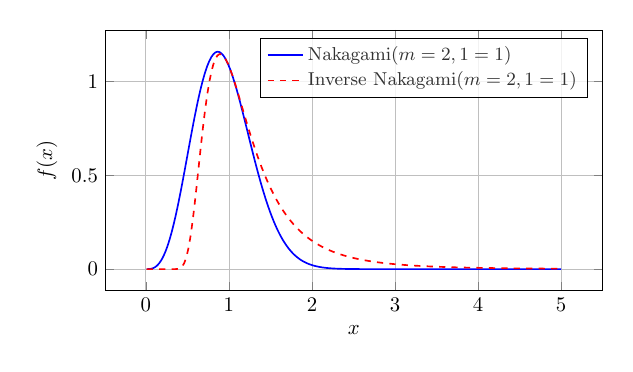}
    \caption{The probability density functions (PDF) for Nakagami and inverse Nakagami distributions. }
    \label{fig:nakagami}
\end{figure}

In cross-entropy (CE/ICE) methods, the use of light-tailed biased densities \(q(\mathbf{u};\mathbf{v})\) — such as Gaussian or vMFNM — can further exacerbate this issue when \(P_{\mathcal{F}}\) is very small. The slow convergence in such cases not only increases computational costs but may also introduce numerical instability during estimation due to the weight degeneracy. To overcome these limitations, we propose an adaptive mixture model that incorporates a heavy-tailed component to accelerate convergence and enhance stability.
  Consider the mixture proposal
\[
q_{\mathrm{safe}}(\mathbf u;\mathbf \phi)
=\lambda\,q(\mathbf u;\mathbf{v})
+(1-\lambda)\,q_H(\mathbf u;\mathbf{\theta}),
\]
where \(q, q_H\) are the corresponding light-tailed and heavy-tailed distributions with corresponding parameters $\mathbf{v}, \mathbf{\theta}$, and \(\lambda\in[0,1]\) controlling the proportion of the light-tailed component.  Drawing samples from \(q_{\mathrm{safe}}\) ensures a more stable and reliable exploration than using \(q\) alone. However, fully utilizing the heavy-tailed kernel will greatly reduce the effective sample size. Therefore, by gradually increasing \(\lambda\) over ICE iterations, the sampler starts with a heavy-tailed component for broad coverage and smoothly transitions toward the optimal light-tailed importance-sampling density, achieving both faster and more robust convergence. For example, when $q$ is a mixture of Gaussian distributions, $q_{H}$ could possibly be the Student-t \cite{ho2012maximum} distribution. Here we consider $q$ to be the vMFNM model in \cite{papaioannou2019improved}, which has the polar form,
\begin{equation}
    \label{vMFNM_model}
    q([r,a];\mathbf{v}) = \sum_{k=1}^{K}\pi_{k}q_{\text{vMFN}}([r, a];\mathbf{v}_{k}),
\end{equation}
where $\pi_{k}\geq 0$ are the weights such that $\sum_{k=1}^{K}\pi_{k} = 1$ and $q_{\text{vMFN}}([r, a];\mathbf{v}_{k}) = q_{\mathcal{N}}(r;m_{k}, \Omega_{k})q_{\text{vMF}}(a;\mu_{k}, \kappa_{k})$. For example, the Nakagami distribution $q_{\mathcal{N}}$ reads the form
\begin{equation}
q_\mathcal{N}\bigl(r; m, \Omega\bigr)
= \frac{2\,m^m}{\Gamma(m)\,\Omega^m}\;r^{2m-1}
\exp\!\Bigl(-\frac{m\,r^2}{\Omega}\Bigr),
\end{equation}
where $\Gamma(\cdot)$ is the Gamma function. As another example, the von-Mises-Fisher distribution has the form
\begin{equation}
q_{\mathrm{vMF}}(\mathbf a;\,\mu,\kappa)
= C_d(\kappa)\,\exp\bigl(\kappa\,\mu^T a\bigr),
\end{equation}
with the normalizing constant
\begin{equation}
C_n(\kappa)
= \frac{\kappa^{\tfrac d2-1}}
       {(2\pi)^{\tfrac d2}I_{\tfrac d2-1}(\kappa)},
\end{equation}
where \(I_{\nu}\) denotes the modified Bessel function of the first kind of order \(\nu\), $d$ is the dimension, and $\|\mu\| = 1$. To construct the heavy-tailed counterpart, we adopt a similar polar decomposition approach. Specifically, we model the radial part with the inverse Nakagami (IN) distribution \cite{louzada2018inverse}, which is a right-heavy-tailed distribution as depicted in Fig. \ref{fig:nakagami} and maintain the angular part with the von-Mises-Fisher distribution. The probability density function of the IN distribution is 
\begin{equation}
    \label{probability_density_function}
    q_{\mathcal{IN}}(r;m, \Omega) = \frac{2\,m^{m}}{\Gamma(m)\,\Omega^{m}}
\,r^{-(2m+1)}
\exp\!\Bigl(-\frac{m}{\Omega\,r^{2}}\Bigr),
\quad x>0.
\end{equation}
Therefore, the novel heavy-tailed distribution $q_{H}$ can be written as 
\begin{equation}
\label{inverse_nakagami}
q_{H}([r, a];\mathbf{\theta}) = \sum_{k=1}^{K}\pi_{k}q_{\mathcal{IN}}(r;m_{k}, \Omega_{k})q_{\text{vMF}}(a;\mu_{k}, \kappa_{k}),
\end{equation}
where $\pi_{k}$ are exactly the weights in Eq. \eqref{vMFNM_model}. As a consequence, the new mixture model is written as 
\begin{equation}
\label{safe_distribution}
    q_{\text{safe}}([r, a];\mathbf{\phi}) = \sum_{k=1}^{K}\left(\lambda q_{\mathcal{N}}(r;m_{k}^{\mathcal{N}}, \Omega_{k}^{\mathcal{N}}) + (1 - \lambda) q_{\mathcal{IN}}(r;m_{k}^{\mathcal{IN}}, \Omega_{k}^{\mathcal{IN}})\right)q_{\text{vMF}}(a;\mu_{k}, \kappa_{k}),
\end{equation}
where $\mathbf{v} = \{\pi_{k}, m_{k}^{\mathcal{N}}, \Omega_{k}^{\mathcal{N}}\}_{k=1}^{K}$ and $\mathbf{\theta} = \{\pi_{k}, m_{k}^{\mathcal{IN}}, \Omega_{k}^{\mathcal{IN}}\}_{k=1}^{K}$.
Then, the stochastic optimization problem \eqref{intermediate_optimization_problems} becomes 
\begin{equation}
\label{novel_optimization}
    \hat{\mathbf{v}}_{t} = \arg\max_{\mathbf{v}\in V}\sum_{i=1}^{N}\ln q(\mathbf{u}_{i};\mathbf{v})W_{t}(\mathbf{u}_{i}, \hat{\mathbf{\phi}}_{t-1}),
\end{equation}
where $\{\mathbf{u}_{i}\}_{i=1}^{N}\sim q_{\text{safe}}(\cdot;\hat{\mathbf{\phi}}_{t-1})$ and $W_{t}(\mathbf{u}_{i}, \widehat{\mathbf{\phi}}_{t-1}) = \frac{p_{t}(\mathbf{u}_{i})}{q_{\mathrm{safe}}(\mathbf{u}_{i};\widehat{\mathbf{\phi}}_{t-1})} $ . That is to say, at each iteration, we will generate samples from the new mixture model $q_{\text{safe}}$ and then solve the corresponding optimization problem. In practice, the parameters $\mathbf{v}$ are determined using the EM algorithm with closed-form updates. However, the parameters of the inverse Nakagami distribution $\{m^{\mathcal{IN}}, \Omega^{\mathcal{IN}}\}$ are not determined, and we have to manually set them to consistently generate high-quality samples to explore the tail space. Note that the mean and mode points for Nakagami and inverse Nakagami distributions \cite{yacoub1999higher, louzada2018inverse} are 
\begin{equation}
\begin{split}
&r_{\mathrm{mode}}^{\mathcal{N}}
= \sqrt{\frac{2m^{\mathcal{N}}-1}{2m^{\mathcal{N}}}\,\Omega^{\mathcal{N}}}, 
\quad (m^{\mathcal{N}}>\tfrac12),\quad
r_{\mathrm{mean}}^{\mathcal{N}}
= \frac{\Gamma\bigl(m^{\mathcal{N}}+\tfrac12\bigr)}{\Gamma(m^{\mathcal{N}})}\,
  \sqrt{\frac{\Omega^{\mathcal{N}}}{m^{\mathcal{N}}}}.\\
&r_{\mathrm{mode}}^{\mathcal{IN}}
= \sqrt{\frac{2m^{\mathcal{IN}}}{2m^{\mathcal{IN}}+1}\,\frac{1}{\Omega^{\mathcal{IN}}}}, 
\quad (m^{\mathcal{IN}}>\tfrac12),\quad 
r_{\mathrm{mean}}^{\mathcal{IN}}
= \frac{\Gamma\bigl(m^{\mathcal{IN}}-\tfrac12\bigr)}{\Gamma(m^{\mathcal{IN}})}\,
  \sqrt{\frac{m^{\mathcal{IN}}}{\Omega^{\mathcal{IN}}}},
  \end{split}
\end{equation}
A natural choice to select $m^{\mathcal{IN}}, \Omega^{\mathcal{IN}}$ is to match the modes between these two distributions, i.e., $r_{\mathrm{mode}}^{\mathcal{N}} = r_{\mathrm{mode}}^{\mathcal{IN}}$. However, when $m\to \frac{1}{2}^{+}$, $r_{\mathrm{mode}}\to 0$. Therefore, to avoid numerical instability, we instead choose $r_{\mathrm{mean}}^{\mathcal{N}} = r_{\mathrm{mode}}^{\mathcal{IN}}$ as $r_{\mathrm{mean}}^{\mathcal{N}} \approx r_{\mathrm{mode}}^{\mathcal{N}}$ when $m^{\mathcal{N}}$ is large, leading to the choice of $\Omega^{\mathcal{IN}}$ as 
\begin{equation}
\label{inverse_gamma_update}
    \Omega^{\mathcal{IN}} = \frac{2m^{\mathcal{IN}}}{2m^{\mathcal{IN}}+1}\left(\frac{\Gamma(m^{\mathcal{N}})}{\Gamma(m^{\mathcal{N}} + \frac{1}{2})}\right)^{2}\frac{m^{\mathcal{N}}}{\Omega^{\mathcal{N}}}.
\end{equation}
Note that the tail thickness of the inverse Nakagami distribution $q_{\mathcal{IN}}\sim x^{-2m^{\mathcal{IN}}-1}$. Therefore,  we can fix a comparatively small $m^{\mathcal{IN}}$ to ensure a heavy tail, leading to a larger exploration space and thus accelerating convergence. However, when the dimension increases, the prior measure will concentrate on a shell. Therefore, choosing $m^{\mathcal{IN}}$ too small in such cases may decrease the effective sample size. To strike a balance between the extrapolation and exploitation ability, we choose $m^{\mathcal{IN}} = \lceil\sqrt{d}\rceil$. Then, we can use Eq.\eqref{inverse_gamma_update} to determine $\Omega^{\mathcal{IN}}$ at each iteration. 

The remaining problem is to choose the proportion $\lambda$ of $q_{H}$ adaptively to ensure a smooth transition from the heavy-tailed dominated case to the light-tailed dominated case. This is achieved by applying a cosine annealing strategy based on the smoothing parameter $\sigma$. In detail, let 
\begin{equation}
\label{cosine_annealing}
    \lambda = \left\{
    \begin{array}{l}
    0,\, \sigma > M.\\ 
    \frac{1}{2}(1 + \cos(\frac{\pi\sigma}{M})), \sigma \leq M.
    \end{array}
    \right.
\end{equation}

Fig.~\ref{fig:cosine} illustrates the above function for the case $M = 2$. When $\sigma$ is large, all samples are drawn from the heavy-tailed distribution, thereby enlarging the exploration space to cover the failure region. As $\sigma$ decreases to $\sigma \leq M$ during the iteration, the majority of samples are gradually generated from the light-tailed distribution, better approximating the optimal bias density. For the selection of $M$, we choose to set $M = \sigma_{0}$, which means $\lambda_{0} = 0$ and all initial samples are coming from the heavy-tailed kernel to explore more space. This adaptive switching strategy improves numerical stability and requires fewer iterations, as the support of $q_{\text{safe}}$ can more effectively cover that of the optimal bias density, particularly when the failure region is located far from the origin.

For the stop criteria in Eq. \eqref{stop_criteria}, we still choose to use the indicator based on the samples from the vMFNM distribution. When $\lambda=0$, meaning no samples come from vMFNM, we will set $\delta_{W}^{*}$ to infinity to continue the iteration. Overall, our new method, termed Safe-ICE-vMFNM, can be summarized as follows:

\begin{figure}[t]
    \centering
    \includegraphics[width=0.5\linewidth]{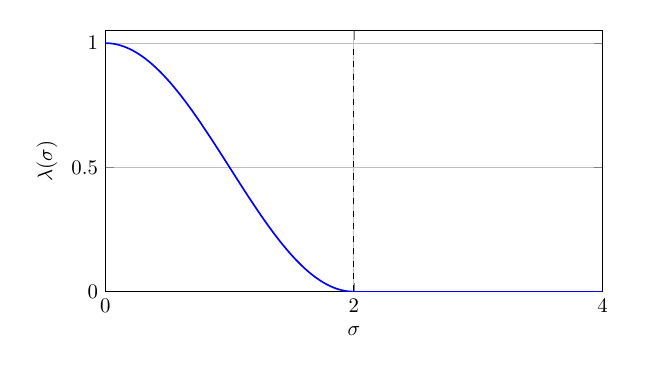}
    \caption{The cosine annealing function when $M = 2$.}
    \label{fig:cosine}
\end{figure}

\begin{itemize}
    \item[1.] Initialize $\sigma_{0}$ and set $\lambda_{0}$ based on \eqref{cosine_annealing}. Choose $K_{0}$ and initialize $\widehat{\mathbf{\phi}}_{0} = \{\widehat{\mathbf{v}}_{0}, \widehat{\mathbf{\theta}}_{0}\}$ based on \eqref{inverse_gamma_update}. Set $t =0$.
    \item[2.] Generate a set of samples $\{\mathbf{u}_{i}\}_{i=1}^{N}$ from $q_{\mathrm{safe}}(\cdot;\widehat{\mathbf{\phi}}_{t})$ and calculate the corresponding LSF values $\{g(\mathbf{u}_{i})\}_{i=1}^{N}$.
    \item[3.] Calculate the CV of weights $\delta_{W_{t}^{*}}$. If it is larger than $\delta^{*}$, go to step 4. Otherwise, go to step 7.
    \item[4.] Calculate $\sigma_{t}$ based on optimization problem \eqref{sigma_optimization}. 
    \item[5.] Set $K(0) = K_{t}, \beta(0) = 1, l(0) = \infty, \widehat{\mathbf{v}}(0) = \widehat{\mathbf{v}}_{t}, j=0$. Using EM to solve stochastic optimization \eqref{novel_optimization}:
    \begin{itemize}
        \item[(a)] Calculate E-step \eqref{expectation}.
        \item[(b)] Set $j = j+1$. Update weights based on \eqref{update_pi}. Update $K(j)$ based on \eqref{update_k}. Update $\beta(j)$ based on \eqref{beta_update}.
        \item[(c)] Update rest parameters based on \eqref{m-step} to obtain $\widehat{\mathbf{v}}(j)$.
        \item[(d)] Calculate $l(j)$ based on \eqref{incomplete_likelihood}.
        \item[(e)] If $|l(j) - l(j-1)| \geq \epsilon \cdot |l(j)|$, no convergence, go to step (a). Otherwise, set $K_{t+1} = K(j), \widehat{\mathbf{v}}_{t+1} = \widehat{\mathbf{v}}(j)$, exit the iteration. 
    \end{itemize}
    \item[6.] Set $t=t+1$. Update $\widehat{\phi}_{t} = \{\widehat{\mathbf{v}}_{t}, \widehat{\theta}_{t}\}$ based on \eqref{inverse_gamma_update}. calculate $\lambda_{t}$ based on \eqref{cosine_annealing}.
    \item[7.] Set $T = t$ and estimate $P_{\mathcal{F}}$ by 
    \begin{equation}
        \widehat{P}_{\mathcal{F}} = \sum_{i=1}^{N}\mathbb{I}_{\Omega_{\mathcal{F}}}(\mathbf{u}_{i})\frac{p(\mathbf{u}_{i})}{q_{\text{safe}}(\mathbf{u}_{i};\widehat{\mathbf{\phi}}_{T})},
    \end{equation}
    where $\{\mathbf{u}_{i}\}_{i=1}^{N}\sim q_{\text{safe}}(\cdot;\widehat{\mathbf{\phi}}_{T})$.
    
\end{itemize}
In the final iteration, all generated samples are used to estimate \(P_{\mathcal{F}}\). This choice does not impact the effective sample size, as our setting ensures that \(\lambda_{T} \to 1\), meaning the sampling distribution in the last iteration coincides with the original light-tailed family. Consequently, the estimator remains statistically valid while fully utilizing the available samples.

\section{Numerical Experiments}
In this section, we evaluate the performance of the proposed Safe-ICE-vMFNM algorithm against the conventional ICE-vMFNM on a series of numerical examples with different orders of magnitude in failure probability.  For each example, we record both the final number of mixture components selected and the total number of iterations $T$ required to compare the performance between the two models. Furthermore, we will use $N=1000$ samples in each ICE iteration to estimate the parameters and the failure probability.  To measure the accuracy and reliability of the probability estimator, we define the coefficient of variance (CV) and the relative MSE error as,
\[
\delta[\widehat{P}_{\mathcal{F}}]
\;=\;
\frac{\sqrt{\mathbb{V}\bigl[\widehat{P}_{\mathcal{F}}\bigr]}}
     {\mathbb{E}\bigl[\widehat{P}_{\mathcal{F}}\bigr]}, \quad \epsilon[\widehat{P}_{\mathcal{F}}] = \frac{\left|P^{MC}_{\mathcal{F}} - \mathbb{E}\bigl[\widehat{P}_{\mathcal{F}}\bigr]\right|}{P_{\mathcal{F}}^{MC}},
\]
where $\widehat{P}_{\mathcal{F}}$ and $P_{\mathcal{F}}^{MC}$ are the importance sampling estimators and MC estimators. To calculate these metrics, we will run each experiment 50 times and calculate the corresponding statistics. Unless otherwise stated, we use a stopping criterion of \(\delta^{*}=1.5\) with a target coefficient of variation \(\delta_{\mathrm{target}}=4\). The initial number of mixtures is set to \(K_{0}=20\). $M$ is specified as $\delta_{0}$. Moreover, the maximum number of iterations for both the outer ICE iterations and the inner EM iterations is set to 20. The prior density $p(x)$ is the standard Gaussian distribution for all experiments.

\subsection{Series system reliability problem}
As an illustrative case, we consider the reliability of a four-component series system. Its limit-state function (LSF) in the two-dimensional standard normal space is given by: 
\begin{equation*}
    g_2(\mathbf{u})
= g(\mathbf{u}) + z, 
\end{equation*}
where 
\begin{equation}
   g(\mathbf{u}) =  \min \left\{
\begin{array}{l}
0.1\,(u_1 - u_2)^2 - \tfrac{1}{\sqrt{2}}(u_1 + u_2) + 3,\;\\ 
0.1\,(u_1 - u_2)^2 + \tfrac{1}{\sqrt{2}}(u_1 + u_2) + 3,\;\\
u_1 - u_2 + \tfrac{7}{\sqrt{2}},\;\\
u_2 - u_1 + \tfrac{7}{\sqrt{2}}
\end{array}
\right\},
\end{equation}
and $z$ is a predefined constant used to control the magnitudes of the failure probabilities. The failure region is defined by $\{g_{2}(\mathbf{u})\leq 0\}$, which has four modes and can be used to test the performance of our adaptive chosen procedure. The ground truth is estimated by MC simulation with $10^{13}$ samples. 
\begin{figure}
    \centering
    \includegraphics[width=0.5\linewidth]{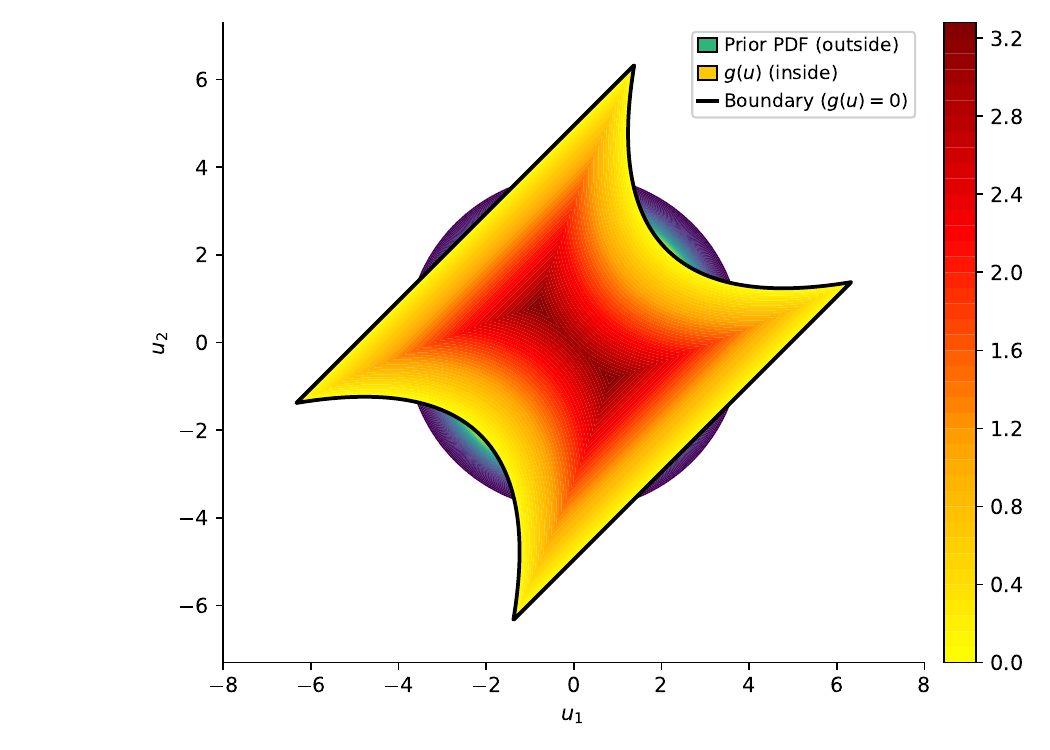}
    \caption{Four-mode problem. The optimal importance density truncated by the failure boundary when $z = 0$.}
    \label{fig:failure_region}
\end{figure}

\begin{figure}
    \centering
    \begin{overpic}[width = 0.45\textwidth]{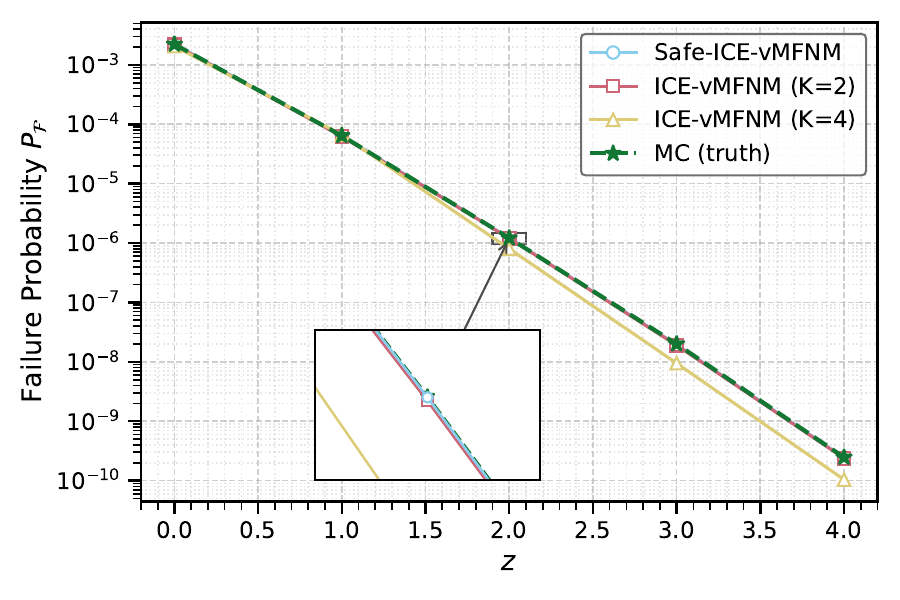}
    \end{overpic}
    \begin{overpic}[width = 0.45\textwidth]{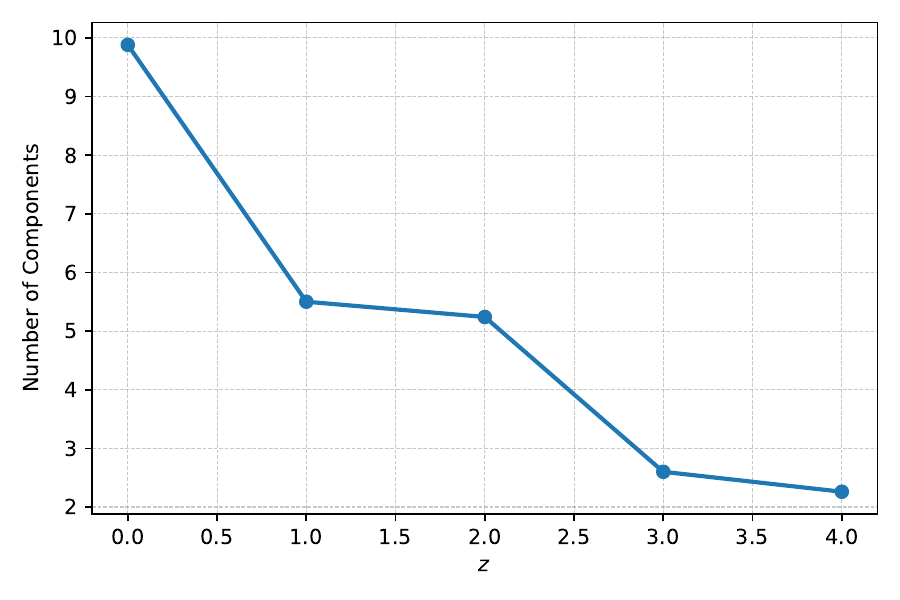}
    \end{overpic}
    \caption{Four-mode problem. Left: the failure probabilities obtained by three methods and the MC estimator for different $z$. Right: The number of components in Safe-ICE-vFMNM for different $z$. All results are reported as the mean of 50 independent runs.}
    \label{fig:pf_4mode}
\end{figure}

To evaluate convergence, we vary the threshold \(z\) from 0 to 4 so as to span a wide range of failure‐probability magnitudes. For ICE-vMFNM, two modes already capture most of the mixture weight as depicted in Fig. \ref{fig:failure_region}, so we test with \(K=2\) and \(K=4\) to assess the impact of mixture size. At each \(z\), we repeat the experiment 50 times and take the sample mean as our approximation of \(P_{\mathcal{F}}\). Figure~\ref{fig:pf_4mode} shows the resulting estimates for all three methods, together with the number of mixtures automatically selected by our approach. While every method yields accurate estimates of \(P_{\mathcal{F}}\), our method attains the highest precision at \(z=4\) and adaptively selects the right number of mixtures as \(z\) grows. 

To obtain a more comprehensive comparison among the three methods, we listed all the statistics in Table \ref{tab:ice-comparison}. For ICE-vMFNM, simply increasing the number of mixtures could lead to worse performance, possibly due to the weight imbalance especially when the failure probability is small. Choosing $K$ to small while covering the majority of mass could lead to better performance, which demonstrates that choosing appropriate number of mixtures is important. However, this drawback is alleviated by our methods, where we adaptively erase the redundant components and thus increase the accuracy and stability, as shown in the table. Combining the mixture proposal with a heavy-tailed kernel, our method achieves faster convergence by requiring fewer iterations, even for very small failure probabilities. This is possibly because the samples generated by the heavy-tailed component can explore a more tailed space and therefore enable a faster decrease of the smoothing parameter. Due to these advantages, our method is not only far more accurate and stable than ICE-vMFNM, but also greatly reduces the computational cost during the iteration process.
\begin{table}[htbp]
  \centering
  \caption{Four-mode problem. The summary of performance for different methods when $z$ changes from 0 to 5.}
  \label{tab:ice-comparison}
  \sisetup{
  table-format           = 1.3,
  table-number-alignment = center,
  table-text-alignment   = center
}
  \begin{tabular}{c *{3}{c c c}}
    \toprule
      & \multicolumn{3}{c}{Safe-ICE-vMFNM}
      & \multicolumn{3}{c}{ICE-vMFNM (K=2)}
      & \multicolumn{3}{c}{ICE-vMFNM (K=4)} \\
    \cmidrule(lr){2-4} \cmidrule(lr){5-7} \cmidrule(lr){8-10}
      {$z$}
      & {$\epsilon[\hat P_f]$} 
      & {$\delta[\hat P_f]$}  
      & {$T[\hat P_f]$}  
      & {$\epsilon[\hat P_f]$} 
      & {$\delta[\hat P_f]$}  
      & {$T[\hat P_f]$}  
      & {$\epsilon[\hat P_f]$} 
      & {$\delta[\hat P_f]$}  
      & {$T[\hat P_f]$}   \\
    \midrule
    0  & \textbf{0.001} & \textbf{0.041} & \textbf{1.0} & 0.019 & 0.084 & 2.3 & 0.071 & 0.122 & 2.1 \\
    1  & \textbf{0.001} & \textbf{0.053} & \textbf{2.0} & 0.029 & 0.111 & 3.6 & 0.035 & 1.175 & 3.5 \\
    2  & \textbf{0.003} & \textbf{0.064} & \textbf{2.2} & 0.013 & 0.369 & 5.0 & 0.341 & 0.888 & 7.1 \\
    3  & \textbf{0.025} & \textbf{0.123} & \textbf{3.1} & 0.049 & 0.186 & 6.1 & 0.526 & 0.915 & 10.8 \\
    4  & \textbf{0.015} & \textbf{0.091} & \textbf{3.3} & 0.037 & 0.117 & 7.5 & 0.578 & 1.003 & 11.7 \\
    \bottomrule
  \end{tabular}
\end{table}

\subsection{Two-dimensional example with three failure regions}
\begin{figure}
    \centering
    \begin{overpic}[width = 0.45\textwidth]{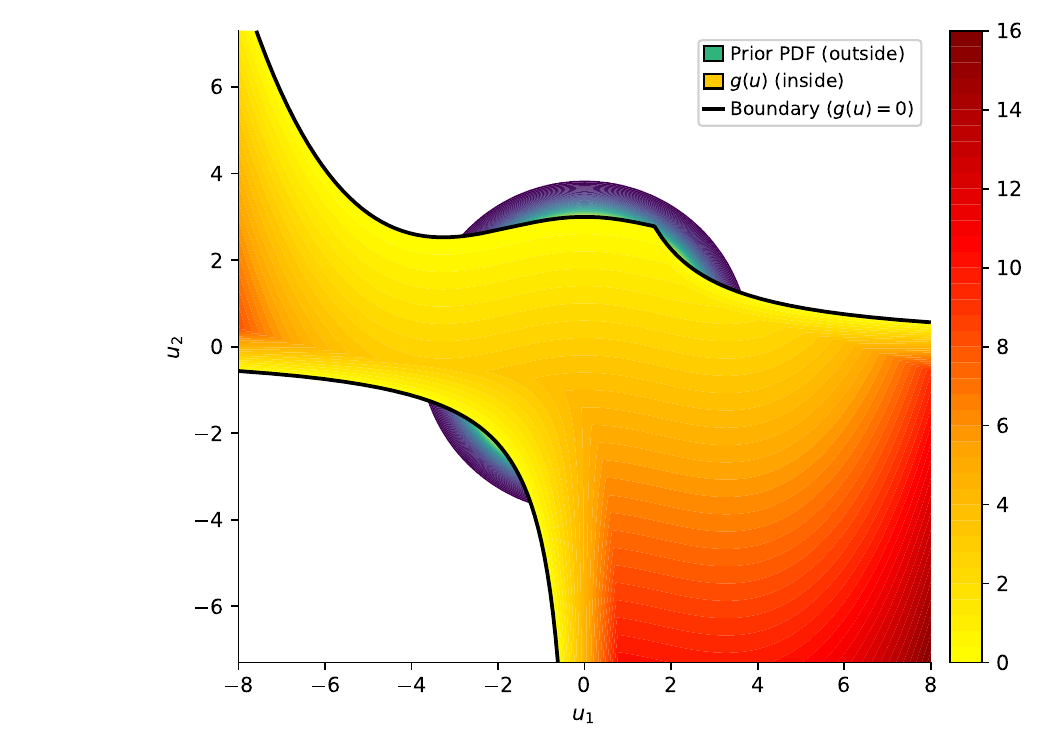}
    \end{overpic}
    \begin{overpic}[width = 0.45\textwidth]{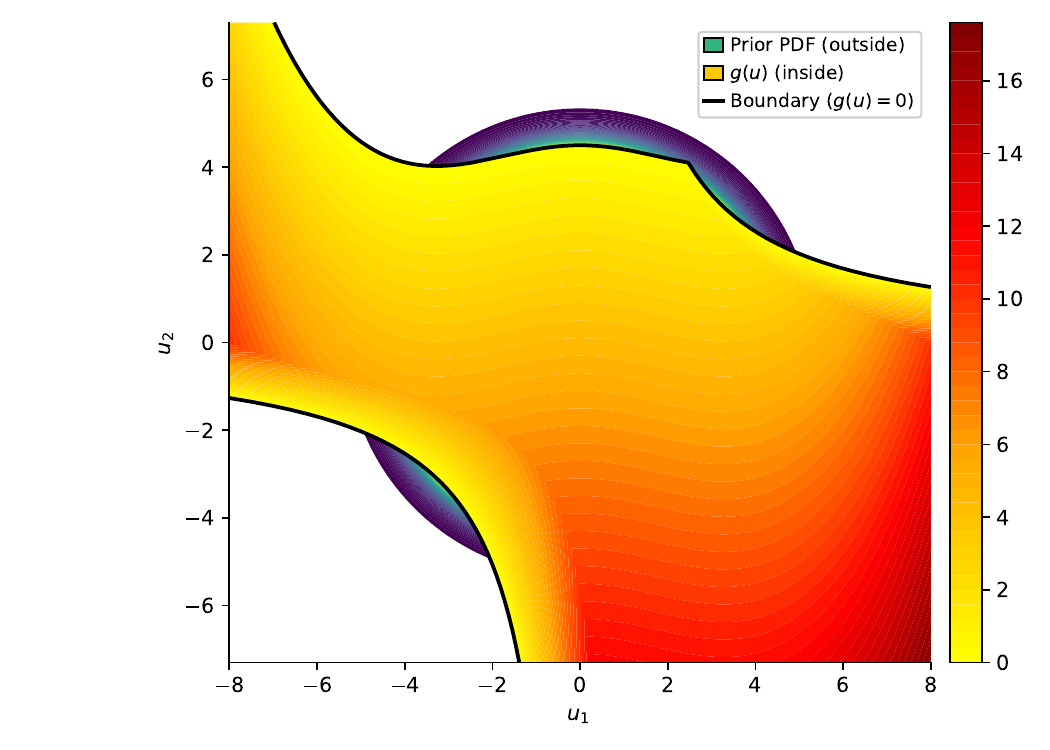}
    \end{overpic}
    \caption{Three-mode problem. The limit-state function and the truncated prior density when $z = 3$ (left) and $4.5$ (right) respectively.}
    \label{fig:failure_region_3mode}
\end{figure}

This example is used to investigate whether our method can identify all the failure regions and thus provide a good estimate of the failure probability. When $d = 2$, the performance function is defined as 
\begin{equation}
    g(\mathbf{u}) = \min\!\left\{\begin{array}{l}
  z - 1 - u_2 + \exp\!\left(-\frac{u_1^{2}}{10}\right) + \left(\frac{u_1}{5}\right)^{4},\\
  \;\frac{z^{2}}{2} - u_1 u_2,
  \end{array}
\right\}
\end{equation}
where $z$ is a prespecified threshold. Here we take $z$ from $3$ to $4.5$ to produce different magnitudes of failure probabilities. The ground truth is obtained by using the MC simulation with $10^{8}$ samples. Note that this example poses significant challenges for adaptive importance sampling procedures. As $z$ increases, the limit-state function undergoes abrupt changes near the third mode in the failure region, as shown in Fig.~\ref{fig:failure_region_3mode}. This makes it easy to overlook the third mode during the iteration process, leading to an insufficient allocation of probability mass and, consequently, a biased failure probability estimate. The two techniques proposed in our method, however, can effectively alleviate such issues, as demonstrated in the following numerical results.

\begin{figure}
    \centering
    \begin{overpic}[width = 0.45\textwidth]{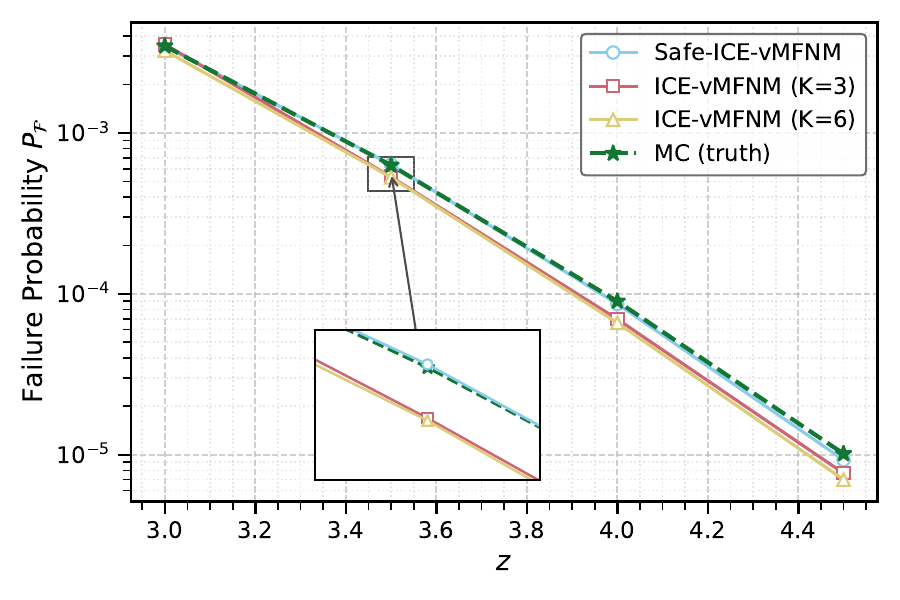}
    \end{overpic}
    \begin{overpic}[width = 0.45\textwidth]{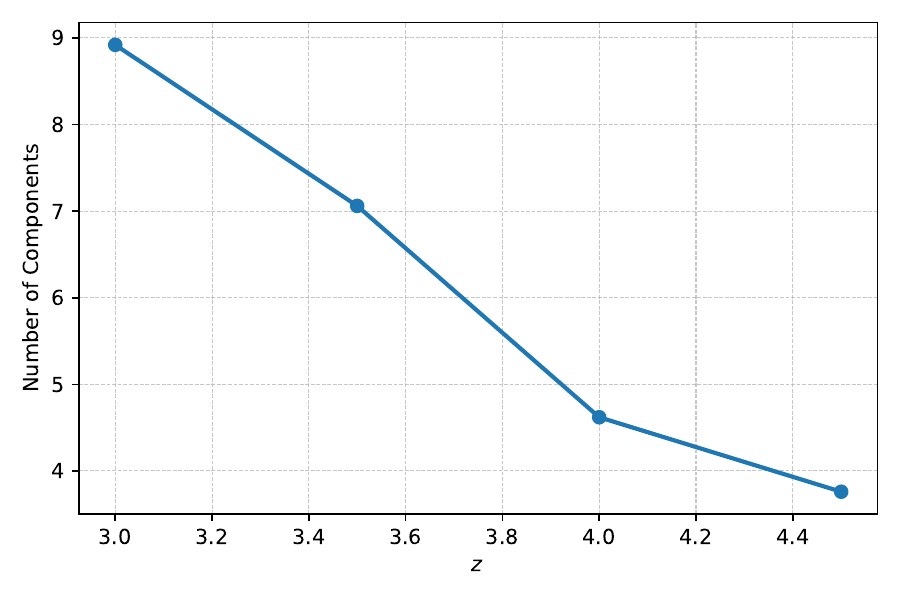}
    \end{overpic}
    \caption{Three-mode problem. Left: the failure probabilities obtained by three methods and the MC estimator for different $z$. Right: The number of components in Safe-ICE-vFMNM for different $z$. All results are reported as the mean of 50 independent runs.}
    \label{fig:pf_3mode}
\end{figure}

For comparison, we initialize ICE-vMFNM with $K=3$ and $6$ to examine whether increasing the number of components improves performance over our adaptive approach. Figure~\ref{fig:pf_3mode} presents the mean failure probability estimates obtained by the three methods, along with the mean number of components selected by our method. The results show that simply increasing $K$ in ICE-vMFNM does not yield better performance. This is because samples drawn from the light-tailed vMFNM model fail to identify the third mode during the iterations, resulting in a clear estimation bias, especially for larger $z$. In contrast, our approach employs a heavy-tailed kernel to explore a broader tail region, increasing the possibility of detecting the third mode. Coupled with adaptive component selection to select enough mixtures, this leads to significantly improved results, as seen in Fig.~\ref{fig:pf_3mode}. 
\begin{table}
  \centering
  \caption{Three-mode problem. Performance summary for different methods.}
  \label{tab:three_mode}
  \sisetup{
    table-format           = 1.3,
    table-number-alignment = center,
    table-text-alignment   = center
  }
  \begin{tabular}{c *{3}{c c c}}
    \toprule
      & \multicolumn{3}{c}{Safe-ICE-vMFNM}
      & \multicolumn{3}{c}{ICE-vMFNM (K=3)}
      & \multicolumn{3}{c}{ICE-vMFNM (K=6)} \\
    \cmidrule(lr){2-4} \cmidrule(lr){5-7} \cmidrule(lr){8-10}
      {$z$}
      & {$\epsilon[\hat P_f]$} & {$\delta[\hat P_f]$} & {$T[\hat P_f]$}
      & {$\epsilon[\hat P_f]$} & {$\delta[\hat P_f]$} & {$T[\hat P_f]$}
      & {$\epsilon[\hat P_f]$} & {$\delta[\hat P_f]$} & {$T[\hat P_f]$} \\
    \midrule
    3.0 & \textbf{0.003} & \textbf{0.046} & \textbf{1.0}
        & 0.024 & 0.607 & 2.0
        & 0.064 & 0.196 & 2.0 \\
    3.5 & \textbf{0.010} & \textbf{0.058} & \textbf{1.5}
        & 0.151 & 0.324 & 2.8
        & 0.158 & 0.139 & 2.5 \\
    4.0 & \textbf{0.038} & \textbf{0.102} & \textbf{2.0}
        & 0.226 & 0.183 & 3.3
        & 0.266 & 0.194 & 3.4 \\
    4.5 & \textbf{0.084} & 0.152 & \textbf{2.2}
        & 0.242 & \textbf{0.074} & 3.7
        & 0.309 & 0.322 & 5.0 \\
    \bottomrule
  \end{tabular}
\end{table}

For completeness, all results are summarized in Table~\ref{tab:three_mode}. Our method consistently generates more accurate and stable estimates. Moreover, by using a heavy-tailed kernel, our method requires fewer iterations to approach the true failure region, thereby reducing the computational cost. These advantages highlight its potential for handling highly complex systems with superior accuracy, robustness, and efficiency.

\subsection{Nonlinear oscillator}
This example considers a hysteretic single-degree-of-freedom oscillator \cite{papaioannou2019improved,wang2016cross} governed by
\begin{equation}
  m\,\ddot x(t)+c\,\dot x(t)+k\!\left[\alpha\,x(t)+\bigl(1-\alpha\bigr)x_y\,z(t)\right]=f(t),
\end{equation}
where $x(t)$, $\dot x(t)$, and $\ddot x(t)$ denote displacement, velocity, and acceleration; $m=6\times10^{4}\,\mathrm{kg}$ is the mass; $k=5\times10^{6}\,\mathrm{N/m}$ is the linear stiffness; $c=2m\zeta\sqrt{k/m}$ is the viscous damping with damping ratio $\zeta=0.05$; $x_y=0.04\,\mathrm{m}$ is the yielding displacement; and $\alpha\in[0,1]$ partitions the restoring force between the elastic part $x$ and the hysteretic part $z$, with $\alpha=0.1$. The internal hysteretic variable follows the Bouc--Wen law
\begin{equation}
  \dot z(t)=\frac{1}{x_y}\!\left[A\,\dot x(t)-\beta\,|\dot x(t)|\,|z(t)|^{\,n-1}z(t)-\gamma\,\dot x(t)\,|z(t)|^{\,n}\right],
\end{equation}
with $A=1$, $\beta=\gamma=0.5$, and $n=3$. The external load models white-noise ground acceleration via a frequency-domain discretization,
\begin{equation}
  f(t)=-\,m\,\sigma \sum_{i=1}^{d/2}\!\big[\,U_i\cos(\omega_i t)+U_{d/2+i}\sin(\omega_i t)\,\big],
\end{equation}
where $\omega_i=i\,\Delta\omega$ with $\Delta\omega=30\pi/d$ and cut-off frequency $\omega_{\mathrm{cut}}=15\pi$, $\sigma=\sqrt{2S\Delta\omega}$ with white-noise intensity $S=0.005\,\mathrm{m^2/s^3}$, $U_i\sim\mathcal N(0,1)$ are independent standard normal coefficients, and $d$ is the number of frequency components. The limit-state function at $t=8\,\mathrm{s}$ is
\begin{equation}
  g(\mathbf{u})=z-x(8\,\mathrm{s})\,,
\end{equation}
where $z$ is the prescribed displacement threshold defining failure. This example is used to test our method when we have no prior knowledge about the failure region. In this experiment, we take $d=10$ and solve the corresponding ODE with the 4th order Runge Kutta method, where the step size is $\Delta t = 0.01$. To compute the true failure probabilities, we apply the MC simulation with $10^{10}$ samples.

\begin{figure}
    \centering
    \begin{overpic}[width = 0.45\textwidth]{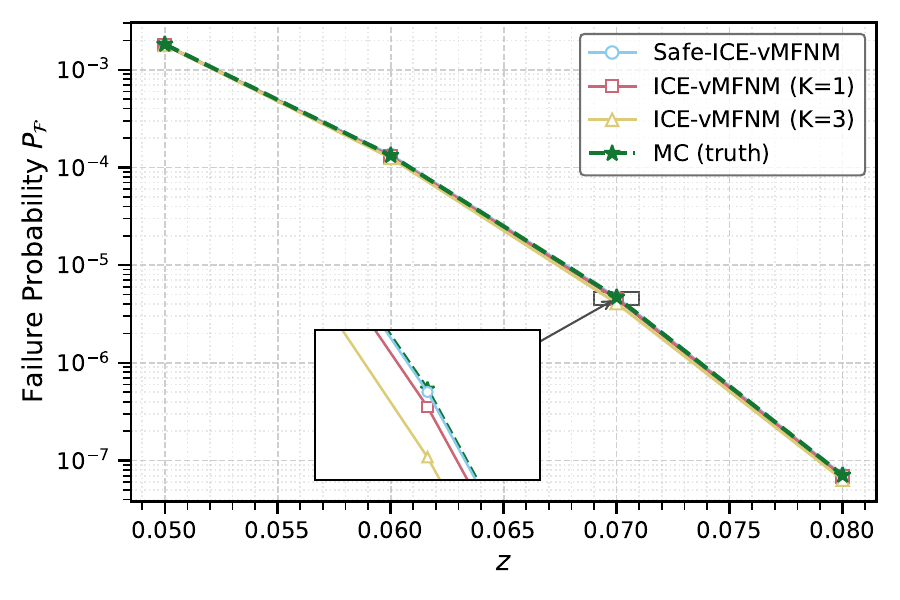}
    \end{overpic}
    \begin{overpic}[width = 0.45\textwidth]{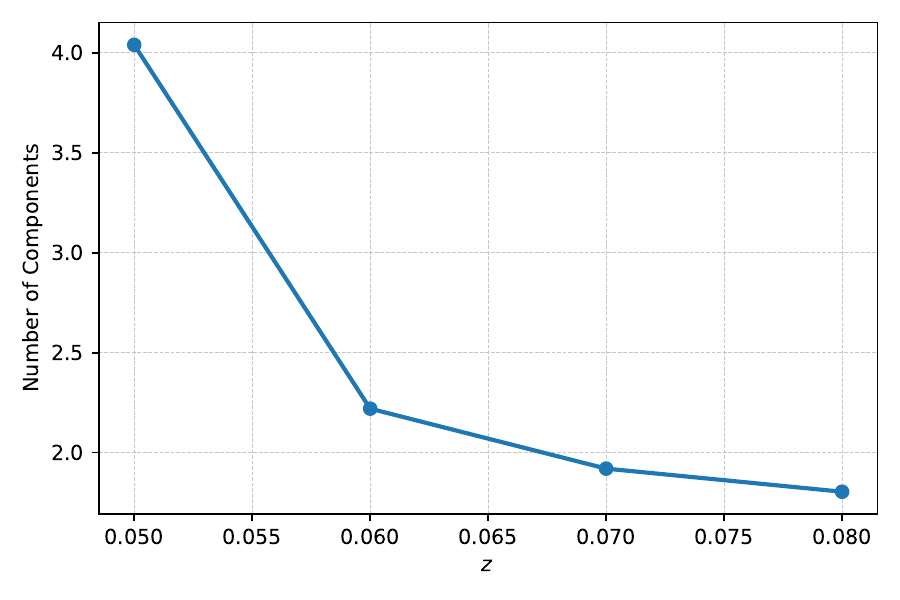}
    \end{overpic}
    \caption{Nonlinear oscillator problem. Left: the failure probabilities obtained by three methods and the MC estimator for different $z$. Right: The number of components in Safe-ICE-vFMNM for different $z$. All results are reported as the mean of 50 independent runs.}
    \label{fig:pf_nonlinear}
\end{figure}

Similarly, we vary the threshold $z$ from $0.05$ to $0.08$ to assess performance across different magnitudes of failure probability. For ICE-vMFNM, we consider mixtures with $K=1$ and $K=3$ to test the impact of different numbers of components for ICE-vMFNM. For each $z$, the experiment is repeated 50 times independently, and the mean estimate is reported; see Fig.~\ref{fig:pf_nonlinear}. Our method automatically selects a sufficient number of components (right panel of Fig.~\ref{fig:pf_nonlinear}), which translates into more accurate probability estimates (left panel of Fig.~\ref{fig:pf_nonlinear}), whereas increasing the number of components in ICE-vMFNM does not consistently improve performance. These advantages are more clearly reflected by Table~\ref{tab:nonlinear}: our estimator remains more accurate even at smaller failure probabilities. In particular, except for the first setting that the results are comparable with each other, its relative error is at least half of that obtained by ICE-vMFNM. Moreover, it attains a smaller coefficient of variation (CV) and converges in fewer iterations, attributable to automatic component reduction and safe exploration via the heavy-tailed proposal. Taken together, these results suggest that the proposed approach handles more complex dynamical systems with greater robustness and lower computational cost.

\begin{table}[htbp]
  \centering
  \caption{Nonlinear oscillator problem. The summary of performance for different methods when $z$ changes from 0.05 to 0.08}
  \label{tab:nonlinear}
  \sisetup{
  table-format           = 1.3,
  table-number-alignment = center,
  table-text-alignment   = center
}
  \begin{tabular}{c *{3}{c c c}}
    \toprule
      & \multicolumn{3}{c}{Safe-ICE-vMFNM}
      & \multicolumn{3}{c}{ICE-vMFNM (K=1)}
      & \multicolumn{3}{c}{ICE-vMFNM (K=3)} \\
    \cmidrule(lr){2-4} \cmidrule(lr){5-7} \cmidrule(lr){8-10}
      {$z$}
      & {$\epsilon[\hat P_f]$} 
      & {$\delta[\hat P_f]$}  
      & {$T[\hat P_f]$}  
      & {$\epsilon[\hat P_f]$} 
      & {$\delta[\hat P_f]$}  
      & {$T[\hat P_f]$}  
      & {$\epsilon[\hat P_f]$} 
      & {$\delta[\hat P_f]$}  
      & {$T[\hat P_f]$}   \\
    \midrule
    0.05  & 0.017 & 0.070 & \textbf{2.0} & 0.016 & 0.082 & 2.3 & \textbf{0.013} & \textbf{0.057} & 2.1 \\
    0.06  & \textbf{0.010} & \textbf{0.091} & \textbf{3.0} & 0.021 & 0.117 & 3.6 & 0.057 & 0.224 & 3.8 \\
    0.07  & \textbf{0.006} & \textbf{0.112} & \textbf{4.1} & 0.036 & 0.153 & 5.4 & 0.130 & 0.312 & 5.7 \\
    0.08  & \textbf{0.00008} & \textbf{0.159} & \textbf{5.5} & 0.028 & 0.178 & 8.1 & 0.099 & 0.397 & 7.2 \\
    \bottomrule
  \end{tabular}
\end{table}

\subsection{System with changing dimension}
In this subsection, we will test a system with two limit-state functions with opposite directions in the normal space, which can be written as 

\begin{equation}
g(\mathbf{u})
= \min\!\Bigl\{
  z - \frac{1}{\sqrt{d}}\sum_{i=1}^d u_i,\;
  z + \frac{1}{\sqrt{d}}\sum_{i=1}^d u_i
\Bigr\},
\end{equation}
where $d$ is the varying dimension and $z$ is user-specified constant.  The failure probability $P_{f}$ is independent of the dimension and can be calculated as $2\Phi(-z)$, where $\Phi(\cdot)$ is the cumulative function of the standard normal distribution. This function is used to test the performance of our method with varying dimensions.

\begin{figure}
    \centering
    \begin{overpic}[width = 0.45\textwidth]{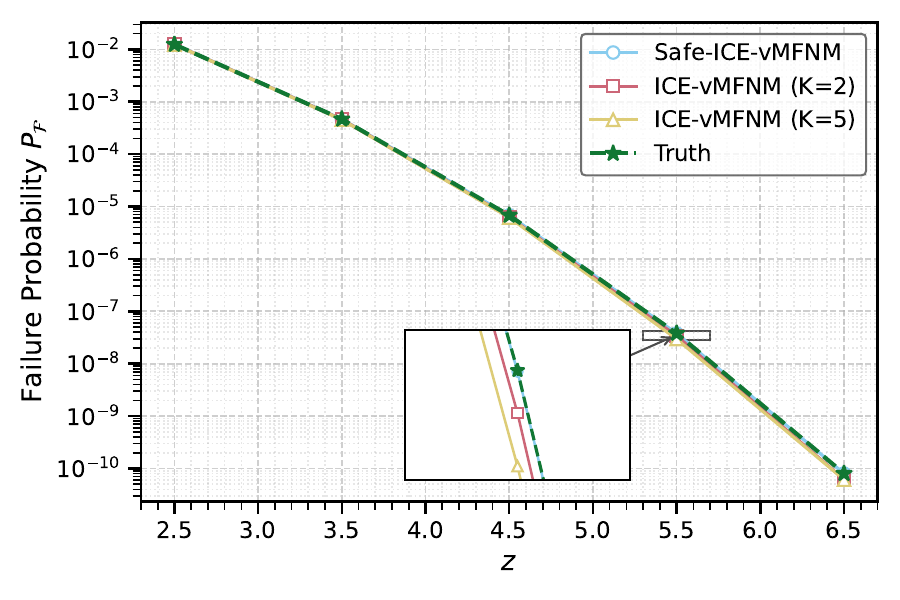}
    \end{overpic}
    \begin{overpic}[width = 0.45\textwidth]{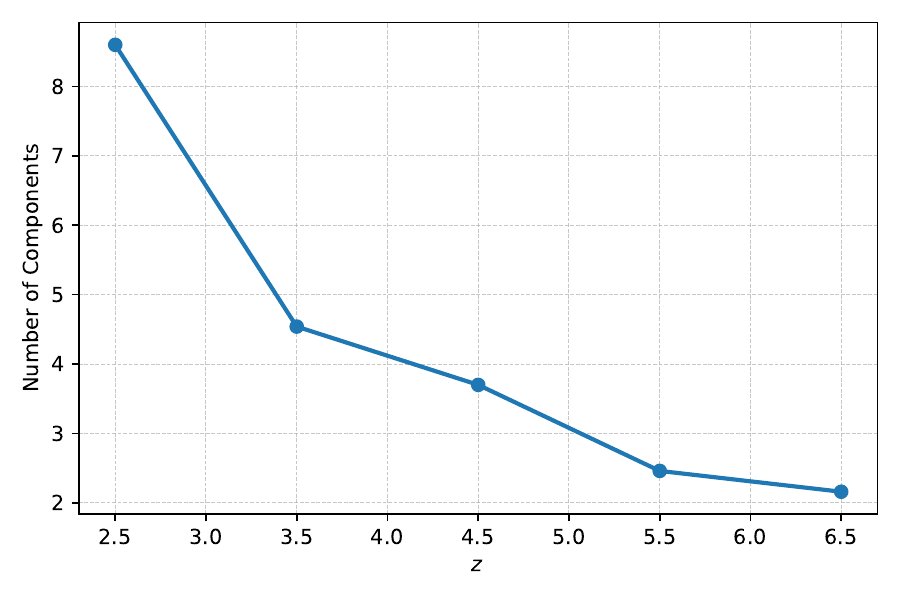}
    \end{overpic}
    \caption{Two-mode problem. Left: the failure probabilities obtained by three methods and the MC estimator for different $z$. Right: The number of components in Safe-ICE-vFMNM for different $z$. All results are reported as the mean of 50 independent runs.}
    \label{fig:pf_2mode}
\end{figure}

We also evaluate \textsc{ICE-vMFNM} with mixture sizes \(K=2\) and \(K=5\) to assess the impact of components across different dimensions. Firstly, for \(d=2\), to test the performance for different magnitudes of failure probabilities, we perform 50 independent trials at each threshold \(z\) from 2.5 to 6.5, compute the sample mean of the estimated failure probability, and plot both the average failure probability estimates and the number of mixture components selected during the iterations in Fig. \ref{fig:pf_2mode}. The conclusion is similar to the previous example. Increasing the number of mixtures of ICE-vMFNM worsens the accuracy as well as the stability. In contrast, our method can adaptively choose the most appropriate number of mixtures to increase the accuracy as depicted in Fig. \ref{fig:pf_2mode}. To comprehensively reveal the phenomenon, we summarize the performance of different methods in Table \ref{tab:two_mode}.

\begin{table}[htbp]
  \centering
  \caption{Two-mode problem. The summary of performance for different methods when $z$ changes from 2.5 to 5.5.}
  \label{tab:two_mode}
  \sisetup{
  table-format           = 1.3,
  table-number-alignment = center,
  table-text-alignment   = center
}
  \begin{tabular}{c *{3}{c c c}}
    \toprule
      & \multicolumn{3}{c}{Safe-ICE-vMFNM}
      & \multicolumn{3}{c}{ICE-vMFNM (K=2)}
      & \multicolumn{3}{c}{ICE-vMFNM (K=5)} \\
    \cmidrule(lr){2-4} \cmidrule(lr){5-7} \cmidrule(lr){8-10}
      {$z$}
      & {$\epsilon[\hat P_f]$} 
      & {$\delta[\hat P_f]$}  
      & {$T[\hat P_f]$}  
      & {$\epsilon[\hat P_f]$} 
      & {$\delta[\hat P_f]$}  
      & {$T[\hat P_f]$}  
      & {$\epsilon[\hat P_f]$} 
      & {$\delta[\hat P_f]$}  
      & {$T[\hat P_f]$}   \\
    \midrule
    2.5  & 0.006 & \textbf{0.037} & \textbf{1.0} & \textbf{0.0008} & 0.052 & 1.5 & 0.004 & 0.053 & 1.2 \\
    3.5  & \textbf{0.001} & \textbf{0.046} & \textbf{1.5} & 0.012 & 0.094 & 2.5 & 0.023 & 0.075 & 2.3 \\
    4.5  & \textbf{0.0003} & \textbf{0.040} & \textbf{2.0} & 0.075 & 0.155 & 3.7 & 0.091 & 0.228 & 3.6 \\
    5.5  & \textbf{0.003} & \textbf{0.047} & \textbf{3.0} & 0.105 & 0.237 & 4.9 & 0.221 & 0.312 & 4.8 \\
    6.5  & \textbf{0.0001} & \textbf{0.055} & \textbf{3.9} & 0.227 & 0.321 & 6.4 & 0.235 & 0.329 & 6.3 \\
    \bottomrule
  \end{tabular}
\end{table}
Table \ref{tab:two_mode} demonstrates that our adaptive ICE-vMFNM consistently outperforms the standard ICE-vMFNM for all tested mixture sizes. By dynamically pruning negligible mixture branches, our method achieves substantially lower relative error and coefficient of variation, yielding both more accurate and more stable failure-probability estimates. Furthermore, the use of a heavy-tailed proposal kernel significantly reduces computational cost—especially when estimating extremely small probabilities—thereby enhancing the robustness of the estimation process. These advantages are particularly valuable for expensive limit-state evaluations, such as those arising from PDE models, where the true number of mixture components is unknown and repeated model solves incur high computational expense.

\begin{figure}[htbp]
    \centering
    \begin{overpic}[width = 0.33\textwidth]{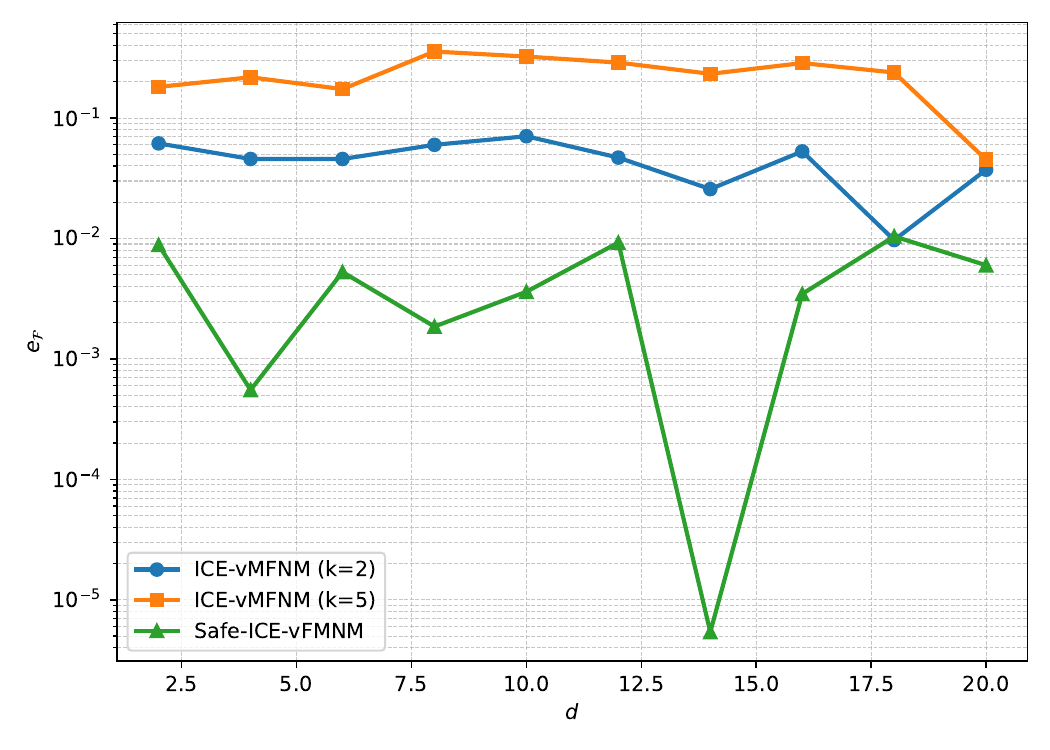}
    \end{overpic}
    \begin{overpic}[width = 0.33\textwidth]{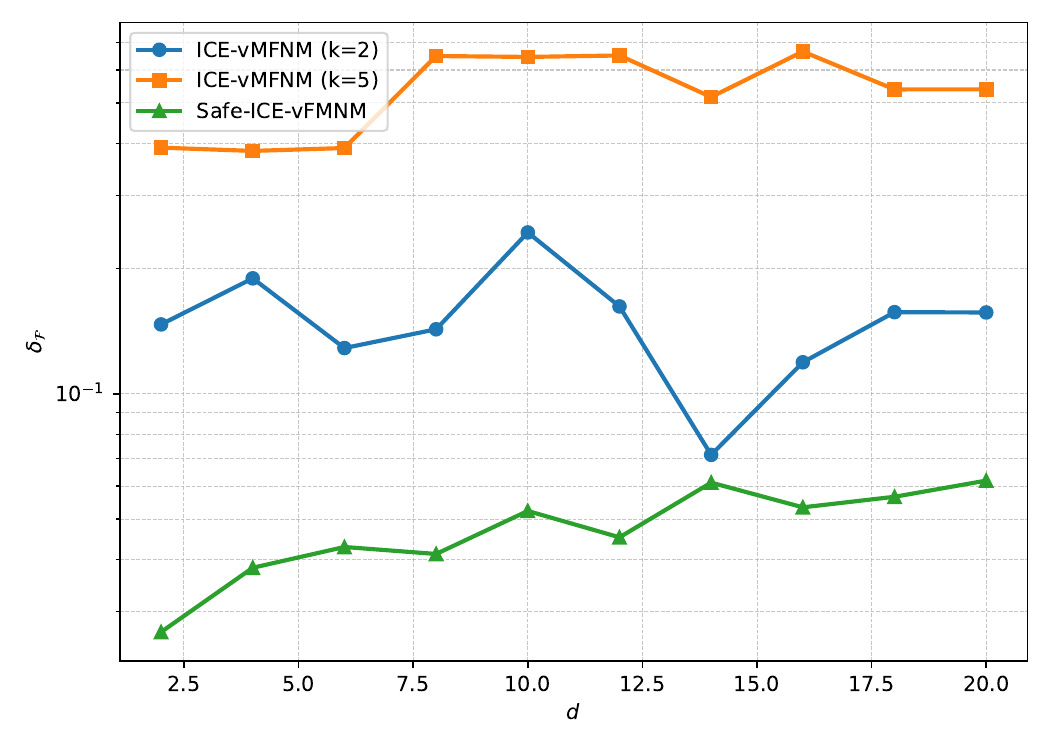}
    \end{overpic}
    \begin{overpic}[width = 0.33\textwidth]{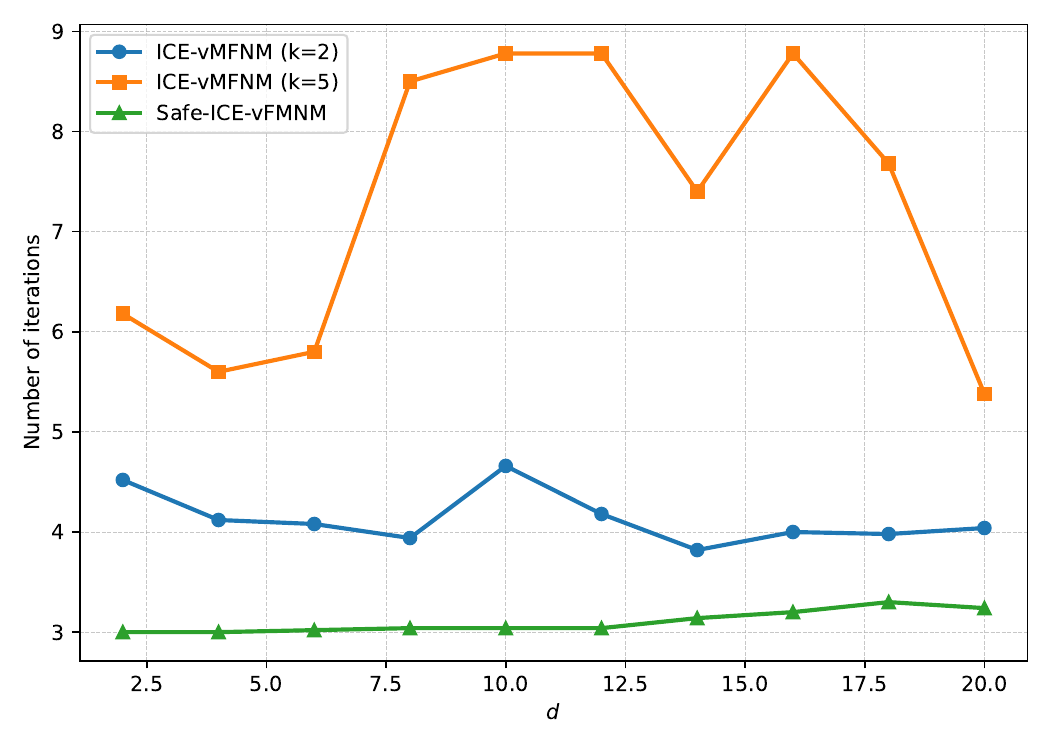}
    \end{overpic}
    \caption{Two-mode problem. Left: the relative error for different methods. Middle: the coefficient of variation for different methods. Right: the number of iterations required for different methods. All results are reported as the mean of 50 independent runs when $d$ increases from 2 to 20. Current $z$ is 5.5.}
    \label{fig:2mode_dimension}
\end{figure}

To isolate the effect of dimensionality, we fix \(z = 5.5\)—which yields a true failure probability of \(3.797\times10^{-8}\)—and sweep the dimension \(d\) from 2 to 20. To stabilize estimates in this mid‐dimensional regime, we increase the per‐trial sample size to 2000, while keeping all other settings as before. Figure~\ref{fig:2mode_dimension} reports, for each method, the relative error, coefficient of variation (CV), and iteration count. Across low to moderate dimensions, our adaptive scheme consistently delivers higher accuracy, smaller CVs, and faster convergence, showing that its performance is largely insensitive to changes in dimension for low to middle dimensional problems. 

\begin{figure}[htbp]
    \centering
    \includegraphics[width=0.4\linewidth]{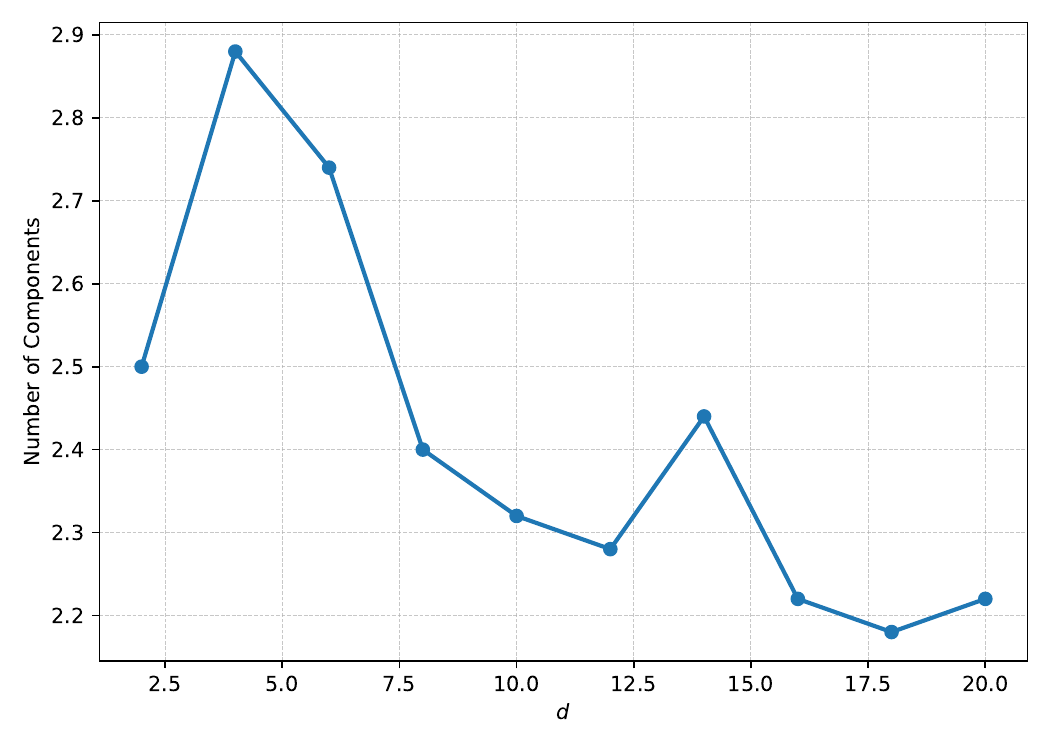}
    \caption{Two-mode problem. The mean number of components selected by our method across different dimensions.}
    \label{fig:dimension_component}
\end{figure}

This advantage stems from two facts. As \(d\) increases further, the Gaussian measure concentrates on a thin shell of radius \(\sqrt{d}\), so even tiny failure regions lie almost exclusively on that shell. In this setting, an overly heavy‐tailed proposal spreads too much mass outside the failure region, causing the effective sample size (ESS) to collapse and sampling efficiency to degrade. Our algorithm avoids this by gradually tempering the tail thickness as \(d\) grows, preserving a high ESS while still covering the critical shell. A secondary but important factor in maintaining efficiency is our dynamic determination of the number of mixture components as depicted in Fig. \ref{fig:dimension_component}, which ensures that only the most relevant modes are included. This adaptive control of both tail weight and component count secures accuracy and stability throughout the low‐to‐moderate dimensional range, where balancing exploration and efficiency matters most. In contrast, simply increasing the number of components in ICE-vMFNM does not incur higher accuracy and stability, which verifies that adaptive selection of components is necessary in our procedure.

Although our method continues to perform well as the dimension $d$ increases, the proposed mixture tends to collapse to the original distribution. In high dimensions, the failure region typically concentrates on a thin shell. Consequently, even for extremely small failure probabilities, most of the probability mass lies on this shell. In such cases, a heavy-tailed radial component offers little benefit. Instead, choosing the original light-tailed distribution with our sparse EM can also work well. Nevertheless, in most practical settings, especially when combined with dimension-reduction techniques \cite{yang2025active,uribe2021cross}, our method can still be applied.

\subsection{Heat transfer problem}

\begin{figure}[htbp]
    \centering
    \begin{overpic}[width = 0.45\textwidth]{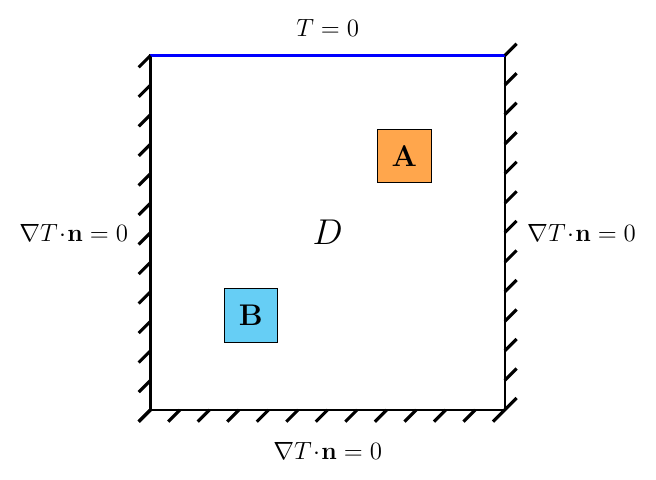}
    \end{overpic}
    \raisebox{.35cm}{
    \begin{overpic}[width=0.29\textwidth]{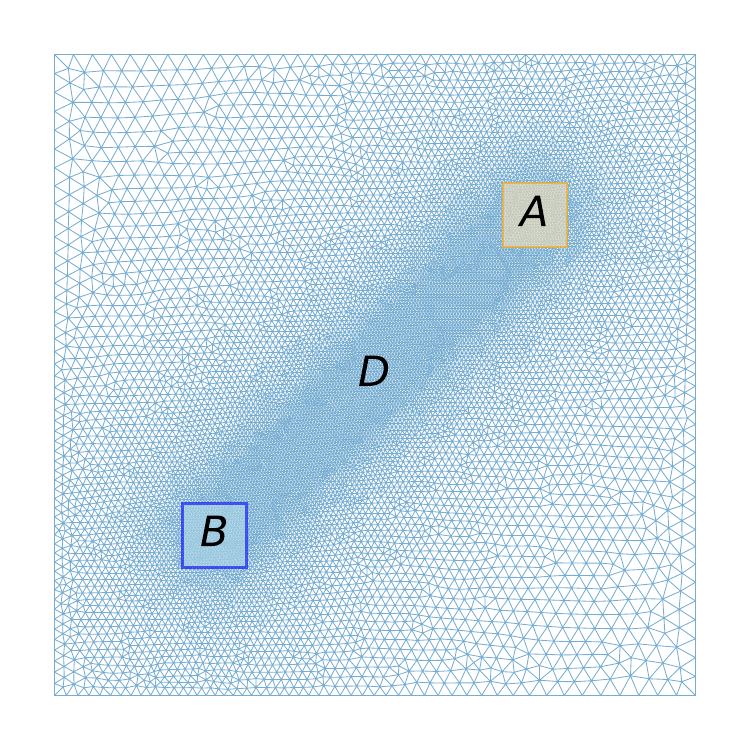}
    \end{overpic}
}
    \caption{Heat conduction: domain and boundary condition (left) and finite element mesh (right).}
    \label{fig:heat_conduction}
\end{figure}

Finally, we test a  heat conduction problem \cite{zhang2025ice} in $D = (-0.5, 0.5)\;m\times (-0.5, 0.5)\;m$; we use finite element simulations. This problem is used to test the performance of our method for real scenarios. Suppose the temperature field $T(\mathbf{x})$ is governed by the following equation:
\begin{equation}
    -\nabla \cdot \left(\kappa(\mathbf{x})\nabla T(\mathbf{x})\right) = I_{A}(\mathbf{x})Q,
\end{equation}
where the boundary condition is divided into two parts: zero Neumann boundary on the top part and zero Dirichlet boundary in the rest of the boundaries, as illustrated in the left figure of Fig.\ref{fig:heat_conduction}. Moreover, we add a heat source term in a localized domain $A = (0.2, 0.3)\;m\times (0.2, 0.3)\; m$ with a constant $Q = 2000 \;W/m^{2}$. The indicator function $I_{A}$ is used to take values in $A$. The permeability field $\kappa(\mathbf{x})$ is a lognormal random field, defined as 
\begin{equation}
    \kappa(\mathbf{x}) = \exp(a_{\kappa} + b_{\kappa} f(\mathbf{x})),
\end{equation}
where $f(\mathbf{x})$ is a standard Gaussian random field with covariance function:
\begin{equation}
    k(\mathbf{x}, \mathbf{x}^{\prime}) = \exp(-\|\mathbf{x} - \mathbf{x}^{\prime}\|^{2}/l^{2}),
\end{equation}
where $l=0.2$ in this paper. The coefficients $a_{\kappa}, b_{\kappa}$ are set such that the mean and the standard deviation of $\kappa(\mathbf{x})$ are $\mu = 1\; W/^{\circ} C \; m$ and $\sigma = 0.3 \; W/^{\circ} C \; m$ respectively. Using the expansion optimal linear estimation \cite{betz2014numerical}(EOLE) method, the random field can be approximated by 
\begin{equation}
    \hat{f}(\mathbf{x}) = \sum_{i=1}^{M} \frac{U_i}{\sqrt{l_i}} \, \boldsymbol{\phi}_i^{\top} \mathbf{C}_{\mathbf{x}\mathbf{\xi}},
\end{equation}
where $U_{i}$ are independent standard normal random variables, $\mathbf{C}_{\mathbf{x}\mathbf{\xi}}$ is a vector with elements $k(\mathbf{x}, \mathbf{\xi}_{i}), i = 1, \ldots, n$, where $\{\mathbf{\xi}_{i},\ldots, \mathbf{\xi}_{n}\}$ are the predefined grid points; $(l_{i}, \phi_{i})$ are the corresponding eigenvalues and eigenvectors of $C_{\mathbf{\xi}\mathbf{\xi}}$ with elements $C^{(k,l)}_{\mathbf{\xi}\mathbf{\xi}} = k(\mathbf{\xi}_{k}, \mathbf{\xi}_{l})$. Generally, it is recommended that the size of the grid is $1/3-1/2$ of $l$ for the exponential covariance function. Hence, we choose the grid size to be $0.1m$, resulting in 121 grid points. For this experiment, we choose the number of terms $M$ to be 10, resulting in over $90\%$ accuracy in the approximation. The limit-state function $g$ is therefore defined as 
\begin{equation}
    g(\mathbf{u}) = 10 - \frac{1}{|B|}\int_{B}T(\mathbf{u})d\mathbf{u},
\end{equation}
where $B$ is a squared domain located in $(-0.3, -0.2)\;m \times (-0.3, 0.2)\; m$. The forward model is solved by the finite element method, with 25040 $T_{3}$ elements, which is displayed in the right figure of Fig.\ref{fig:heat_conduction}. To evaluate the ground truth of the failure probability, we use the subset simulation method to run this experiment 50 times with 2000 samples in each level and calculate the mean as the final approximated ground truth. The reference failure probability is $4.69\times 10^{-7}$.

To investigate the impact of the number of components in ICE-vMFNM, we also set $K = 1$ and $3$, and compare the performance with safe-ICE-vMFNM. All experiments are run 50 times independently. The estimated values of failure probabilities of safe-ICE-vMFNM, ICE-vMFNM (K=1) and ICE-vMFNM (K=3) are $4.72\times 10^{-7}, 4.63\times 10^{-7}$ and $4.26\times 10^{-7}$ respectively. Therefore, our method can even achieve a more accurate estimate in this case, as reported in Table.\ref{tab:heat}.
Moreover, in terms of stability and computational cost, our method achieves a more stable and efficient estimate, with a CV of $0.035$ and a mean iteration number of 3.3. This advantage still stems from the automatic pruning for increasing stability and heavy-kernel-based sampling for enhancing accuracy. Finally, the mean number of components of our model is 1.7, which further explains the good performance of our method. 
  
  \begin{table}[htbp]
  \centering
  \caption{Heat conduction problem. The summary of performance for different methods.}
  \label{tab:heat}
  \sisetup{
  table-format           = 2.5,
  table-number-alignment = center,
  table-text-alignment   = center
}
  \begin{tabular}{l c c c}
    \toprule
    Method & $\epsilon[\hat P_f]$ & $\delta[\hat P_f]$ & $T[\hat P_f]$ \\
    \midrule
    Safe-ICE-vMFNM        & \textbf{0.006}   & \textbf{0.035} & \textbf{3.0} \\
    ICE-vMFNM (K=1)       & 0.013 & 0.044 & 3.3 \\
    ICE-vMFNM (K=3)       & 0.092   & 0.303 & 4.7 \\
    \bottomrule
  \end{tabular}
  
\end{table}

\section{Conclusion}
In this paper, we proposed using a cross-entropy penalized EM algorithm to update the weights in order to erase the redundant components in this mixture. This procedure increases the accuracy as well as the stability. For extremely small failure probabilities, we proposed a new mixture that couples the light-tailed and the heavy-tailed distributions to explore the tail space efficiently. To strike a balance between exploration and exploitation, we adopted the cosine-annealing strategy to adjust these two terms, which finally accelerates the convergence rate and thereby reduces the computational cost.  However, one limitation of our method is that for high dimensions, the heavy-tailed kernel will degrade to the original distribution due to the measure concentration. To test the effectiveness of our method, we used different examples. The results consistently showed that our method is more robust and accurate than the original one for different magnitudes of failure probabilities.

\section*{Acknowledgments}

This work was supported by the MURI grant (FA9550-20-1-0358), the ONR Vannevar Bush Faculty Fellowship (N0001422-1-2795), and the U.S. Department of Energy, Advanced Scientific Computing Research program, under the Scalable, Efficient and Accelerated Causal Reasoning Operators, Graphs and Spikes for Earth and Embedded Systems (SEA-CROGS) project, (DE-SC0023191).
Additional funding was provided by GPU Cluster for Neural PDEs and Neural Operators to support MURI Research and Beyond, under Award \#FA9550-23-1-0671.

\section*{Appendix A}
Here we derive how to update the weights using the \textbf{Q}-function defined in Eq.\eqref{Q-function}. We first write the  Lagrangian $\tilde{J}$ as 
\begin{equation}
\begin{split}
    \tilde{J}(\pi, \mathbf{v}) = \sum_{i=1}^{N}W_{i}\sum_{k=1}^{K}\gamma_{k}^{(i)}\log\pi_{k}q(\mathbf{u}_{i};\mathbf{v}_{k}) + \beta\sum_{i=1}^{N}W_{i}\sum_{k=1}^{K}\pi_{k}\ln \pi_{k} + \lambda(\sum_{k=1}^{K}\pi_{k}-1).
\end{split}
\end{equation}
Take the first derivative of $\tilde{J}$ with respect to $\pi_{k}$ and set it to zero, we can get:
\begin{equation}
\label{derivation}
  \sum_{i=1}^{N}W_{i}\gamma_{k}^{(i)} + \beta \left(\sum_{i=1}^{N}W_{i}\right)\pi_{k}(\ln\pi_{k} + 1) + \lambda \pi_{k} = 0.
\end{equation}
 Summing $k$ leads to 
 \begin{equation}
 \label{lambda}
\lambda = -\sum_{k=1}^{K}\sum_{i=1}^{N}W_{i}\gamma_{k}^{(i)} - \beta \left(\sum_{i=1}^{N}W_{i}\right)\left(1 +\sum_{k=1}^{K}\pi_{k}\ln\pi_{k}\right)
 \end{equation}
Substituting Eq.\eqref{lambda} into Eq.\eqref{derivation}, we have 
\begin{equation}
    \sum_{i=1}^{N}W_{i}\gamma_{k}^{(i)} + \beta \left(\sum_{i=1}^{N}W_{i}\right)\pi_{k}\ln\pi_{k}-\left(\sum_{s=1}^{K}\sum_{i=1}^{N}W_{i}\gamma_{s}^{(i)}\right)\pi_{k} - \beta \left(\sum_{i=1}^{N}W_{i}\right)\left(1 +\sum_{s=1}^{K}\pi_{s}\ln\pi_{s}\right) \pi_{k} = 0.
\end{equation}
Therefore, we have
\begin{equation}
    \pi_{k}^{\text{new}} = \frac{\sum_{i=1}^{N}W_{i}\gamma_{k}^{(i)}}{\sum_{s=1}^{K}\sum_{i=1}^{N}W_{i}\gamma_{s}^{(i)}} + \beta\left(\frac{\sum_{i=1}^{K}W_{i}}{\sum_{s=1}^{K}\sum_{i=1}^{N}W_{i}\gamma_{s}^{(i)}}\right)\pi_{k}^{\text{old}}\left(\ln \pi_{k}^{\text{old}} - \sum_{s=1}^{K}\pi_{s}^{\text{old}}\ln\pi_{s}^{\text{old}}\right),
\end{equation}
where the first term on the right-hand side is exactly the EM update of weights. 
\bibliographystyle{unsrt}  
\bibliography{references}

\end{document}